\documentclass[times,review,10pt]{elsarticle}

\usepackage{amssymb}
\usepackage{amsmath}
\usepackage{array}
\usepackage{caption,color}
\usepackage{subfigure}
\usepackage{float}
\usepackage{textcomp}
\usepackage{amsmath}
\usepackage{stfloats}
\usepackage{url}
\usepackage{verbatim}
\usepackage{amssymb}
\usepackage{bm,amsfonts}
\usepackage{longtable}
\usepackage{textcomp}
\usepackage{braket}
\usepackage{multirow}
\usepackage{cases}
\usepackage{graphicx}
\usepackage{tabularx}
\usepackage{algorithm}
\usepackage{algorithmic}
\usepackage{epstopdf}
\usepackage{setspace}

\newtheorem{defn}{Definition}

\newcommand\bmx{\bm{\mathcal{X}}}
\newcommand\bmy{\bm{\mathcal{Y}}}
\newcommand\bme{\bm{\mathcal{E}}}
\newcommand\bms{\bm{\mathcal{S}}}
\newcommand\bmu{\bm{\mathcal{U}}}
\newcommand\bmv{\bm{\mathcal{V}}}
\newcommand\bmi{\bm{\mathcal{I}}}
\newcommand\wbmi[1]{\widehat{\boldsymbol{#1}}^{(i)}}
\newcommand\Cnnn{\mathbb{C}^{n_1\times n_2\times n_3}}

\journal{Journal of \LaTeX\ Templates}

\begin{document}

\begin{frontmatter}

%% Title, authors and addresses

%% use the tnoteref command within \title for footnotes;
%% use the tnotetext command for theassociated footnote;
%% use the fnref command within \author or \affiliation for footnotes;
%% use the fntext command for theassociated footnote;
%% use the corref command within \author for corresponding author footnotes;
%% use the cortext command for theassociated footnote;
%% use the ead command for the email address,
%% and the form \ead[url] for the home page:
%% \title{Title\tnoteref{label1}}
%% \tnotetext[label1]{}
%% \author{Name\corref{cor1}\fnref{label2}}
%% \ead{email address}
%% \ead[url]{home page}
%% \fntext[label2]{}
%% \cortext[cor1]{}
%% \affiliation{organization={},
%%             addressline={},
%%             city={},
%%             postcode={},
%%             state={},
%%             country={}}
%% \fntext[label3]{}

\title{High-dimensional low-rank three-way tensor autoregressive time
series predictor}

%% use optional labels to link authors explicitly to addresses:
%%
\author{Haoning Wang,\quad  Liping Zhang\corref{mycorrespondingauthor}}
    \cortext[mycorrespondingauthor]{Corresponding author. \\
    \emph{E-mail addresses}: whn22@mails.tsinghua.edu.cn (H. Wang),  lipingzhang@tsinghua.edu.cn (L. Zhang).}
  \address{ Department of Mathematical Sciences, Tsinghua University, Beijing
        100084, China}

%% Abstract
\begin{abstract}
%% Text of abstract
Recently, tensor time-series forecasting has gained increasing attention, whose core requirement is how to perform dimensionality reduction.
In this paper, we establish a least square optimization model by combining tensor singular value decomposition (t-SVD) with autoregression (AR) to forecast third-order tensor time-series, which has great benefit in computational complexity and dimensionality reduction. We divide such an optimization problem using fast Fourier transformation and t-SVD into four decoupled subproblems, whose variables include regressive coefficient, f-diagonal tensor, left and right orthogonal tensors, and propose an efficient forecasting algorithm via alternating minimization strategy, called Low-rank Tensor Autoregressive Predictor (LOTAP), in which each subproblem has a closed-form solution. Numerical experiments indicate that, compared to Tucker-decomposition-based algorithms, LOTAP achieves a speed improvement ranging from $2$ to $6$ times while maintaining accurate forecasting performance in all four baseline tasks. In addition, this algorithm is applicable to a wider range of tensor forecasting tasks because of its more effective dimensionality reduction ability.
\end{abstract}

%% Keywords

\end{frontmatter}

%% Add \usepackage{lineno} before \begin{document} and uncomment
%% following line to enable line numbers
%% \linenumbers

%% main text
%%

%% Use \section commands to start a section
\section{Introduction}

With the exponential growth of data volume, real world time series data often exhibit complex structures that require representation of multiple features or channels \cite{weron2008forecasting}. These complex structures lead to higher-order time series structures.
In the context of time series data, zero to second-order series refer to scalar, vector, and matrix formats of time series data, respectively, which are also known as lower-order time series. In contrast, higher-order time series, also known as tensor-valued time series, encompass third-order and higher-order time series data.

Among these structures, third-order data is the most common form of the higher-order data. Many real-world time-series instances naturally have third-order data structures. For instance, in color video data, each frame consists of RGB images, which inherently have three dimensions. Similarly, in the modeling of meteorological and ocean data, a snapshot of meteorological information at a specific point in time can be effectively represented using three dimensions: longitude, latitude, and variability.

Time series forecasting is a crucial problem in the field of time series data analysis, with numerous applications in various domains such as finance, meteorology, and healthcare. In particular, third-order time series forecasting, has gained significant attention due to the increasing prevalence of multidimensional data in real world scenarios.
A prominent challenge to the third-order time series forecasting problem is the high order nature of the data \cite{cai2015facets}.
The increase in dimensionality can lead to the ``dimensional catastrop", intensifying the computational load of algorithms and complicating data analysis due to the intrinsic complexity of correlations.

To address the challenge associated with the high dimensional nature of the data described above, a technique that can reduce the dimensionality of the higher-order data is required. This necessity naturally leads us to the realm of tensor decomposition techniques. Tensor decomposition is a potent tool for extracting critical information from tensor data. It includes methods such as CANDECOMP/PARAFAC (CP) decomposition \cite{CP2013}, Tucker decomposition \cite{tucker}, and tensor singular value decomposition (t-SVD) \cite{KILMER2011}. These techniques extract latent variables or components that capture the most salient features and eliminate redundancy from the original data. Leveraging these advantages, tensor decomposition-based forecasting methods can manage multiple tensor time series simultaneously and achieve commendable prediction performance.

A significant body of the prior work on higher-order time series forecasting problem is based on tensor Tucker decomposition.
This preference for Tucker decomposition is due to the high computational complexity associated with CP decomposition, which is NP-hard, and the short existence of the t-SVD algorithm.
Moving forward, the recent study by \cite{wang2021high} introduced the LATAR model, which extends the vector autoregressive (AR) model to tensor form by considering the low-rank Tucker decomposition estimation of the transfer tensor. The Multilinear Orthogonal AR (MOAR) and Multilinear Constrained AR (MCAR) models \cite{moar2018jing} were developed for forecasting higher-order time series. These models incorporate a set of projection matrices to obtain the potential core tensor under specific constraints and generalize the traditional AR model to tensor form. This is achieved by combining the Tucker decomposition with the AR model to facilitate higher-order time-series analysis and forecasting. Another approach, presented by \cite{bhtarima2020shi}, is the block Hankel tensor ARIMA (BHT-ARIMA) model, which builds upon the MCAR model by applying the multi-way delay embedding transform (MDT) technique to the tensor time series.

Despite the above mentioned forecasting algorithms based on Tucker decomposition are effective, they have several limitations.
First, these methods assume that the tensor data under Tucker decomposition is low-rank, but as we will show in Section \ref{Section: Numerical}, this may not hold in the real world.
Second, Tucker decomposition necessitates  matrix decomposition after matrixing each mode of the tensor, resulting in significant computational overhead. The core requirement of tensor time-series forecasting is how to perform dimensionality reduction. However, such methods have limitations in terms of computational cost, with iteration complexity of approximately $O(n^3r)$, where $n$ and $r$ are the dimension and rank of original tensor data.

The limitations inherent in the existing studies have provided the impetus to search for an improved model for predicting third-order time series. t-SVD, as a recently proposed tensor decomposition algorithm, its application in various fields has only gradually come into our view in recent years. We observed a surge in the development of tensor completion algorithms utilizing t-SVD and its bootstrapped tensor nuclear norm in recent years, particularly in the context of solving third-order tensor completion problems \cite{zhang2016exact,tensorfactorization2018zhou,zhang2019nonlocal,gilman2020online,li2021t,zhang2021multiscale}. This tensor decomposition algorithm, derived directly from the tensor-tensor product (t-product) \cite{kilmer2011factorization}, has exhibited comparable performance to the conventional Tucker decomposition-based tensor completion algorithm \cite{filipovic2015tucker,zhao2015bayesian,liu2019low} and has demonstrated superiority in terms of filling problems and broader applicability.

\begin{figure}[!ht]
    \centering
    \includegraphics[width=0.9\textwidth]{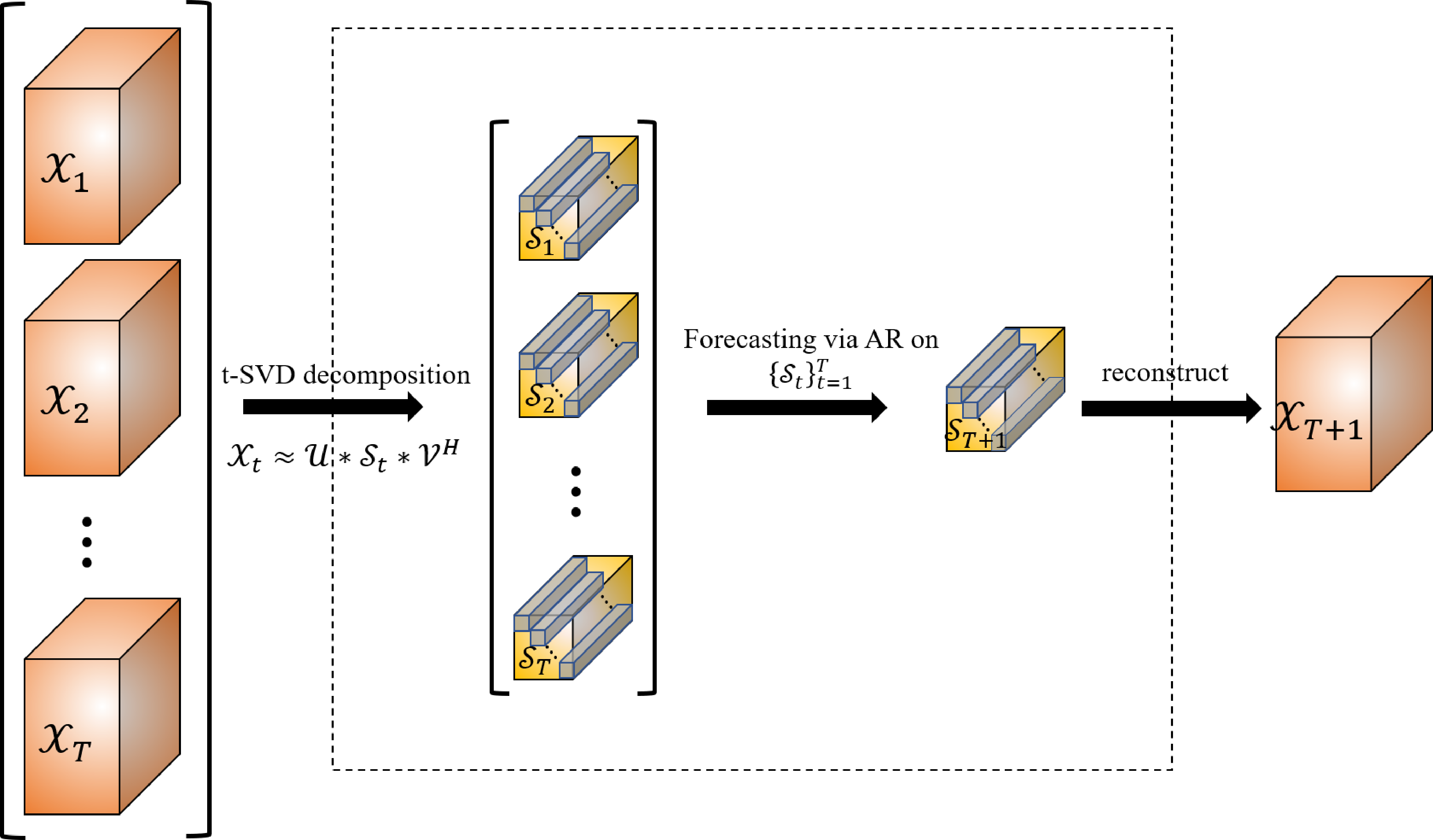}
    \caption{\small{Schematic illustration of the core idea of our proposed LOTAP method.}}
    \label{fig:flow chart of LOTAP}
\end{figure}

Given the promising performance of t-SVD in tensor completion, we utilize the truncated t-SVD into third-order time series forecasting problem for the first time, and then propose an efficient algorithm, \textbf{Lo}w-rank \textbf{T}ensor \textbf{A}utoregressive \textbf{P}redictor (LOTAP), which solves a least squares optimization problem formulated by AR and truncated t-SVD.
Unlike the Tucker-decomposition-based algorithms, our proposed LOTAP can always preserve the intrinsic low-rank structure, thus exhibiting a wider range of applications and leading to a greater advantage in computational complexity.
Moreover, by transforming the original data into the Fourier domain, LOTAP can circumvent the need to directly compute the t-product differential in the original optimization problem. Thereby, we propose a new alternating minimization algorithm to solve LOTAP, in which each subproblem has a closed-form solution.
% Compared with other state-of-the-art algorithms, our proposed LOTAP algorithm significantly reduces the computational overhead and exhibits a wider range of applications.
Experimental results on synthetic and real datasets demonstrate that our LOTAP model and proposed solution method outperform a few state-of-the-art methods. The fundamental concept of our LOTAP model is depicted in Figure \ref{fig:flow chart of LOTAP}.

The core contributions of this paper are summarised as follows.
\begin{enumerate}
    \item We utilize truncated t-SVD to capture the critical information about the third-order time-series data and then  build an efficient forecasting algorithm, LOTAP. To the best of our knowledge, this is the first time to introduce truncated t-SVD together with AR model into the field of tensor time-series forecasting problem.
    \item LOTAP formulates third-order tensor time-series forecasting problem  as a least squares optimization and divides such an optimization problem using fast Fourier transformation and t-SVD into four decoupled subproblems, whose variables include regressive coefficient, f-diagonal tensor, left and right orthogonal tensors. A new alternating minimization algorithm is proposed  to solve such problem, in which each subproblem has a closed-form solution. LOTAP not only approximately reduce the iteration complexity from $O(n^3r)$ to $O(n^3+n^2r^2)$ but also has a broader application compared to forecasting algorithms based on Tucker decomposition due to its maintaining the nature of the low-rank structure of tensor time-series data.

\end{enumerate}

The rest of the paper is structured as follows. In Section II, we briefly recall previous work related to  higher-order time-series forecasting based on AR model. In Section III, some basic notation of tensor algebra and preliminaries are provided. The proposed LOTAP model and alternative minimizing algorithm are described in Section IV. Numerical experiments on synthetic and real-world data are reported in Section V, and the conclusions are given in Section VI.

\section{Preliminaries}
In this section, we first introduce the basic tensor notation with tensor algebra and then review the AR model. Except for some specific cases, we represent a scalar with ordinary letters, e.g., $x$; a column vector with lowercase bold letters, e.g., $\boldsymbol{x}$; a matrix with uppercase bold letters, e.g., $\boldsymbol{X}$; and a tensor with Euler script letters, e.g., $\bmx$.

\paragraph{Tensor Algebra}
An N-order tensor is a multi-linear structure in $\mathbb{C}^{n_1\times \dots\times n_N}$. For a third-order tensor $\bmx\in\Cnnn$, we use the Matlab notation $\bmx(k,:,:),\bmx(:,k,:)$ and $\bmx(:,:,k)$ to denote the $k$-th horizontal, lateral, and frontal slices, $\bmx(:,i,j),\bmx(i,:,j)$ and $\bmx(i,j,:)$ to denote the mode-1, mode-2 and mode-3 fibers. Specifically, the frontal slice $\bmx(:,:,k)$ is also denoted as $\boldsymbol{X}^{(k)}$. The Frobenius norm of $\bmx$ is defined as $$\lVert\bmx\rVert_F := \sqrt{\sum_{i,j,k}\lvert\bmx(i,j,k)\rvert^2}.$$

For any tensor $\bmx\in\Cnnn$, we denote $\widehat{\bmx}\in \Cnnn$ as the results of the fast Fourier transform (FFT) along the third dimension, {i.e.}, $\widehat{\bmx}:= \mathcal{H}(\bmx) = \text{fft}(\bmx,[],3)$ in the Matlab command. We can compute $\bmx$ from $\widehat{\bmx}$ by inverse fast Fourier transform (IFFT) along the third dimension, {i.e.}, $\bmx =\mathcal{H}^{-1}(\widehat{\bmx}) =\text{ifft}(\widehat{\bmx},[],3)$. Let $\text{conj}(\bmx)$ denote the complex conjugate of $\bmx$. Let $\overline{\boldsymbol{X}}$ denote the block diagonal matrix of the tensor, where the $k$-th diagonal block of $\overline{\boldsymbol{X}}$ is the $k$-th frontal slice of $\widehat{\bmx}$, {i.e.},
$$
    \overline{\boldsymbol{X}}:=\operatorname{bdiag}(\widehat{\mathcal{X}})=\left[\begin{array}{llll}
    \widehat{\boldsymbol{X}}^{(1)} & & & \\
    & \widehat{\boldsymbol{X}}^{(2)} & & \\
    & & \ddots & \\
    & & & \widehat{\boldsymbol{X}}^{\left(n_3\right)}
    \end{array}\right].
$$
It is shown in \cite{TRPCA2019Lu} that
$$\lVert\bmx\rVert_F = \frac{1}{\sqrt{n_3}}\lVert\overline{\boldsymbol{X}}\rVert_F=\frac{1}{\sqrt{n_3}}\lVert\widehat{\bmx}\rVert_F.$$

Now we define three matrix block-based operators, and then lead to the definition of the t-product between two third-order tensor and the definition of truncated t-SVD.

Let $\text{bcirc}(\bmx)$ and $\text{bvec}(\bmx)$ denote the block circular matrix and the block vectorization matrix, respectively,  constructed by $\boldsymbol{X}^{(k)}$, {i.e.}
$$
    \operatorname{bcirc}(\bmx):=\left[\begin{array}{cccc}
    \boldsymbol{X}^{(1)} & \boldsymbol{X}^{\left(n_3\right)} & \cdots & \boldsymbol{X}^{(2)} \\
    \boldsymbol{X}^{(2)} & \boldsymbol{X}^{(1)} & \cdots & \boldsymbol{X}^{(3)} \\
    \vdots & \vdots & \ddots & \vdots \\
    \boldsymbol{X}^{\left(n_3\right)} & \boldsymbol{X}^{\left(n_3-1\right)} & \cdots & \boldsymbol{X}^{(1)}
    \end{array}\right],
$$
$$
    \operatorname{bvec}(\bmx):=\left[\begin{array}{c}
    \boldsymbol{X}^{(1)}\\
    \boldsymbol{X}^{(2)}\\
    \vdots\\
    \boldsymbol{X}^{\left(n_3\right)}
    \end{array}\right].
$$
It is shown in \cite{TRPCA2019Lu} that the block circular matrix can be block diagonalized by FFT, {i.e.},
$$
(\bm{F}_{n_3} \otimes \bm{I}_{n_1}) \cdot \text{bcirc}(\bmx) \cdot (\bm{F}_{n_3}^{-1} \otimes \bm{I}_{n_2})=\overline{\boldsymbol{X}},
$$
where $\bm{F}_n$ is the $n\times n$ discrete Fourier matrix, $\bm{I}_n$ is the $n\times n$ identity matrix and $\otimes$ denotes the Kronecker product.

The inverse operator $\text{bvfold}(\cdot)$ takes $\text{bvec}(\cdot)$ to a tensor form: $\text{bvfold}(\text{bvec}(\bmx)) = \bmx$. Then the t-product can be defined as follows:
\begin{defn} [\textbf{t-product} \cite{KILMER2011}]
    Let $\bmx\in\Cnnn$ and $\bmy\in\mathbb{C}^{n_2\times n_4\times n_3}$. Then the t-product $\bmx *\bmy$ is defined to be a tensor of size $n_1\times n_4\times n_3$,
    \begin{equation}\label{tproducdef}
        \bmx *\bmy = \mathrm{bvfold}(\mathrm{bcirc}(\bmx)\cdot\mathrm{bvec}(\bmy)).
    \end{equation}
\end{defn}
Clearly, we have the following equivalence:
$$\bmx *\bmy = \bm{\mathcal{Z}}\Leftrightarrow \overline{\boldsymbol{X}}\cdot\overline{\boldsymbol{Y}}=\overline{\boldsymbol{Z}}.$$

\begin{defn} [\textbf{conjugate transpose} \cite{TRPCA2019Lu}]
    The conjugate transpose of a tensor $\bmx\in\Cnnn$ is the tensor $\bmx^H\in\Cnnn$ obtained by conjugate transposing each of the frontal slices and then reversing the order of transposed frontal slices 2 through $n_3$.
\end{defn}

\begin{defn} [\textbf{f-diagonal tensor} \cite{KILMER2011}]
    A tensor $\bmx$ is called f-diagonal if each frontal slice $\boldsymbol{X}^{(k)}$ is a diagonal matrix.
\end{defn}

\begin{defn} [\textbf{identity tensor} \cite{KILMER2011}]
    The identity tensor $\bmi_{n_1,n_3}\in\mathbb{C}^{n_1\times n_1\times n_3}$ is a tensor whose first frontal slice is an identity matrix and the remaining frontal slices are all zeros.
\end{defn}
It is easy to show that $\bmx *\bmi = \bmi *\bmx = \bmx$ and $\overline{\bmi}$ is a tensor with each frontal slice being identity matrix.
\begin{defn} [\textbf{column-orthogonal tensor}]
A tensor $\bmx\in\Cnnn$ is a column-orthogonal tensor if and only if
$\bmx^H *\bmx = \bmi_{n_2,n_3}$.
\end{defn}

The column-orthogonal tensor is a degeneration of orthogonal tensor \cite{KILMER2011}. It is obvious that a column-orthogonal tensor $\bmx\in\Cnnn$  is an \textbf{orthogonal tensor} if and only if $n_1=n_2$.

\begin{defn} [\textbf{tensor singular value decomposition: t-SVD} \cite{KILMER2011}]
    For $\bmx\in\Cnnn$, the t-SVD of $\bmx$ is given by $\bmx = \bmu *\bms *\bmv^H$, where $\bmu\in\mathbb{C}^{n_1\times n_1\times n_3}$ and $\bmv\in\mathbb{C}^{n_2\times n_2\times n_3}$ are orthogonal tensors, and $\bms\in\Cnnn$ is an f-diagonal tensor, respectively.
\end{defn}
\begin{defn} [\textbf{tensor tubal rank} \cite{tensorfactorization2018zhou}]
    For $\bmx\in\Cnnn$, the tensor tubal rank, denoted as $\text{rank}_t(\bmx)$, is defined by $\text{rank}_t(\bmx) = \max_{1\leq i\leq n_3}\{\text{rank}(\widehat{\boldsymbol{X}}^{(i)}) \}.$
\end{defn}

Based on the definition of the column-orthogonal tensor and the tensor tubal rank, we can extend the definition of t-SVD to the truncated form.
\begin{defn} [\textbf{truncated t-SVD}]
    Let $\bmx\in\Cnnn$ be a tensor and $r=\text{rank}_t(\bmx)$ be the tensor tubal rank of $\bmx$, then truncated t-SVD of the tensor $\bmx$, is $\bmx = \bmu *\bms *\bmv^H$, where $\bmu\in\mathbb{C}^{n_1\times r\times n_3}$ and $\bmv\in\mathbb{C}^{n_2\times r\times n_3}$ are column-orthogonal tensors, and $\bms\in\mathbb{C}^{r\times r\times n_3}$ is an f-diagonal tensor, respectively.
\end{defn}

\begin{figure}[!ht]
    \centering
    \includegraphics[width=0.7\textwidth]{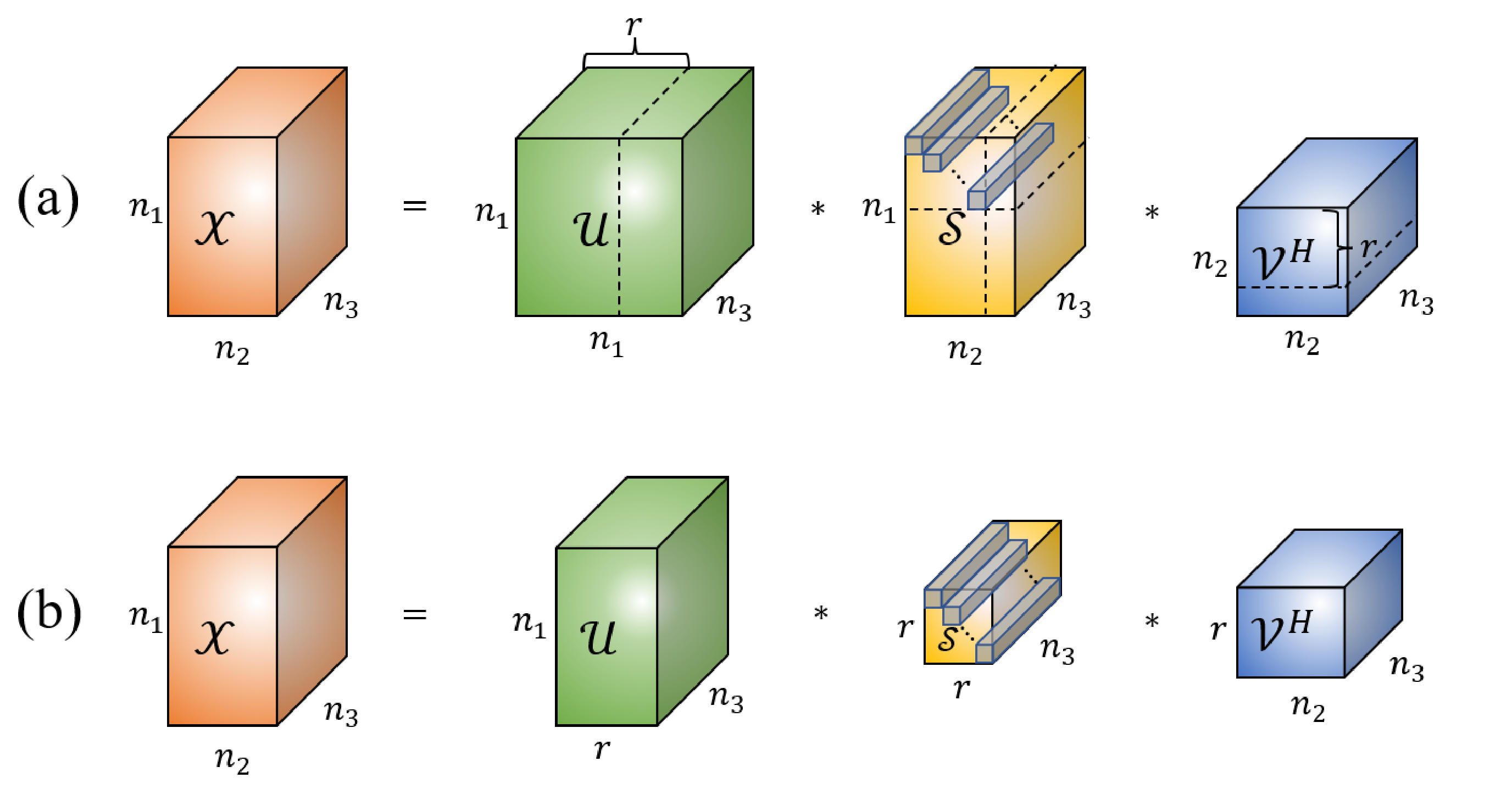}
    \caption{\small{a) t-SVD of tensor $\bmx$; b) truncated t-SVD of tensor $\bmx$ for the truncated rank $r$.}}
    \label{fig:t-SVD and truncated t-SVD}
\end{figure}

Figure \ref{fig:t-SVD and truncated t-SVD} shows the t-SVD and its truncated version. 

{It is worth noting that, besides the standard (deterministic) truncated t-SVD \cite{ahmadi2024adaptive}, randomized variants have also been developed to substantially reduce computational cost \cite{zhang2018randomized,ahmadi2024randomized,chen2023parallel}. Specifically, while the classical t-SVD requires $O(n_1 n_2 n_3 \min\{n_1,n_2\})$ operations, randomized t-SVD reduces the complexity to $O(n_1 n_2 n_3 k)$, where $k$ is a target rank parameter with $k \ll \min\{n_1,n_2\}$. Such randomized approaches thus offer a more efficient alternative for large-scale applications. }

{We stress that in our work, t-SVD is introduced solely for modeling purposes within the optimization formulation; no explicit t-SVD computations are carried out during the actual algorithmic iterations.}

\paragraph{Autoregressive Model}
Let $x_t$ be the actual data value at the time point $t$. Then, $\text{AR}(p)$ model would regard $x_t$  as the linear combination of the past $p$ values, {i.e.},
\begin{equation}
    x_t = \sum_{i=1}^pa_ix_{t-i}+\epsilon_t,
\end{equation}
where the random error $\{\epsilon_t\}$ of the observation is identically distributed with a mean zero and a constant variance, and $\{a_i\}_{i=1}^p$ are the coefficients of AR. AR models are usually used for single time-series forecasting.

\section{Algorithm Description}\label{Section: Alg-des}

In this section, we present a novel algorithm based on the truncated t-SVD for learning and forecasting third-order tensor time series. The algorithm consists of two steps. In the first step, an appropriate temporal parameter is learned by solving an optimization problem based on an AR model and the truncated t-SVD. In the second step, the tensor data at the next time point is forecasted by using the parameters obtained from the first step.

\subsection{Step 1: Tensor Autoregression with Truncated t-SVD}

Let $\bmx\in \mathbb{C}^{n_1\times n_2\times n_3\times T}$ denote the input time series, where $T$ is the number of observed time points. We define the third-order tensor $\bmx_t\in\Cnnn$ as the slice of $\bmx$ at time point $t$.
Using the truncated t-SVD, we extract the f-diagonal tensor $\bms_t\in\mathbb{C}^{r\times r\times n_3}$ containing core features from $\bmx_t$ by jointly factorizing it into column-orthogonal tensors $\bmu\in\mathbb{C}^{n_1\times r\times n_3}$ and $\bmv\in\mathbb{C}^{n_2\times r\times n_3}$, which can be formulated as:
\begin{equation}\label{X_tsvd}
\bmx_t = \bmu \ast \bms_t\ast \bmv^H.
\end{equation}

Note that in \eqref{X_tsvd} we model the tensor time-series data with shared left and right tensor subspace factor. This is because we assume that for adjacent tensor data on the time scale, their subspace factors hardly change. {This assumption is further validated in Subsection~5.3}. Indeed, in the numerical experiments of Section~\ref{Section: Numerical}, we will set smaller values for $T$ to ensure that the factor subspaces from $\bmx_1$ to $\bmx_T$ do not vary much.

To further explore the intrinsic time-series connections within the original data $\{\bmx_t\}_{t=1}^T$, we utilize the f-diagonal tensors $\{\bms_t\}_{t=1}^T$ to capture the most critical information about the original tensors, and construct an AR model in tensor form of order $p$. Specifically, we represent $\bms_t$ as a linear combination of its $p$ most recent predecessors ${\bms_{t-1},\dots,\bms_{t-p}}$, i.e.,
\begin{align}\label{S_ar}
\bms_t &= \sum_{i=1}^pa_i\bms_{t-i} + \bme_t\nonumber\\
       &\overset{a}{=} \sum_{i=1}^pa_i(\bmu^H\ast\bmx_{t-i}\ast\bmv) + \bme_t\nonumber\\
       &= \bmu^H\ast(\sum_{i=1}^pa_i\bmx_{t-i})\ast\bmv + \bme_t,\quad p+1\leq t\leq T,
\end{align}
where $\{a_i\}_{i=1}^p$ are the AR coefficients and $\{\bme_t\}$ are the random errors of the observation, the equality $a$ holds due to $\bmu^H\ast\bmu=\bmi_{n_1,n_3}$ and $ \bmv^H\ast\bmv=\bmi_{n_2,n_3}$. To minimize the forecast error $\bme_t$ for each time point, we formulate the following  least squares optimization problem based on  \eqref{X_tsvd}:
 \begin{align}\label{optim_problem}
\mathop{\min}_{\{\bms_t\},\{a_i \},\bmu,\bmv} & F\left(\{\bms_t\},\{a_i \},\bmu,\bmv\right)=
\sum_{t=p+1}^{T}\lVert\bms_{t}- \sum_{i=1}^pa_i\bms_{t-i}\rVert^2_F \nonumber\\
&\qquad +\varphi\sum_{t=1}^{T}\lVert\bmx_t-\bmu\ast\bms_t\ast\bmv^H\rVert_F^2 \nonumber\\
\mbox{s.t.}\quad &  \bmu^H\ast\bmu=\bmi_{n_1,n_3},~~\bmv^H\ast\bmv=\bmi_{n_2,n_3},\\
& P_{\Omega}(\bms_t) = \bms_t, ~~\forall 1\leq t\leq T,\nonumber
\end{align}
where $\Omega=\{(i,i,j): 1\leq i\leq r,1\leq j\leq n_3\}$ is the set consisting of the indexes of all diagonal elements in each frontal slice and $P_{\Omega}$ is a linear operator that extracts entries in and fills the entries not in with zeros, which guarantees the f-diagonality of $\bms_t$. The optimization problem \eqref{optim_problem} involves minimizing the squared Frobenius norm of the difference between the f-diagonal tensor $\bms_t$ and its linear combination of the $p$ preceding tensors, subject to the constraint that the joint factor orthogonal tensors $\bmu$ and $\bmv$ are unitary. Additionally, the problem includes a regularization term that enforces the reconstructed tensor data to be close to the observed data in the Fourier domain.

A common approach to minimizing \eqref{optim_problem} is through the \textbf{alternating minimization algorithm}. However, the presence of the t-product ``$\ast$" prevents the direct derivation of closed-form solutions. By transforming the data into the Fourier domain, we can alternatively express the cost function of optimization problem \eqref{optim_problem} in the following form:
 \begin{align}\label{fft_transform}
&\sum_{t=p+1}^{T}\lVert\bms_{t}- \sum_{i=1}^pa_i\bms_{t-i}\rVert^2_F +\varphi\sum_{t=1}^{T}\lVert\bmx_t-\bmu\ast\bms_t\ast\bmv^H\rVert_F^2 \\
=& \sum_{t=p+1}^{T}\frac{1}{n_3}\lVert\overline{\boldsymbol{S}_t}- \sum_{i=1}^pa_i\overline{\boldsymbol{S}_{t-i}}\rVert^2_F
+\frac{\varphi}{n_3}\sum_{t=1}^{T}\lVert\overline{\boldsymbol{X}_t}-\overline{\bmu\ast\bms_t\ast\bmv^H}\rVert_F^2 \nonumber\\
=& \sum_{t=p+1}^{T}\frac{1}{n_3}\lVert\overline{\boldsymbol{S}_t}- \sum_{i=1}^pa_i\overline{\boldsymbol{S}_{t-i}}\rVert^2_F
+\frac{\varphi}{n_3}\sum_{t=1}^{T}\lVert\overline{\boldsymbol{X}_t}-\overline{\mathbf{U}}\cdot\overline{\boldsymbol{S}_t}\cdot\overline{\mathbf{V}}^H\rVert_F^2 \nonumber \\
=& \frac{1}{n_3}\left(\sum_{t=p+1}^{T}\lVert\overline{\boldsymbol{S}_t}- \sum_{i=1}^pa_i\overline{\boldsymbol{S}_{t-i}}\rVert^2_F
+\varphi\sum_{t=1}^{T}\lVert\overline{\mathbf{U}}^H\cdot\overline{\boldsymbol{X}_t}\cdot\overline{\mathbf{V}}-\overline{\boldsymbol{S}_t}\rVert_F^2\right). \nonumber
\end{align}
By applying the Fourier transform presented in Eq. \eqref{fft_transform}, we can transform the original optimization problem \eqref{optim_problem} into a new formulation, denoted as \eqref{fft_optim_problem}, which can be efficiently solved by an alternative minimization algorithm.
 \begin{alignat}{2}\label{fft_optim_problem}
\mathop{\min}_{\{\widehat{\bms}_t\},\{a_i \},\widehat{\bmu},\widehat{\bmv}} &\sum_{t=p+1}^{T}\sum_{i=1}^{n_3}\lVert\wbmi{S}_t- \sum_{i=1}^pa_i\wbmi{S}_{t-i}\rVert^2_F\\
&+\varphi\sum_{t=1}^{T}\sum_{i=1}^{n_3}\lVert\wbmi{X}_t-\wbmi{U}\cdot\wbmi{S}_t\cdot(\wbmi{V})^H\rVert_F^2 \nonumber\\
\mbox{s.t.}\quad &  (\wbmi{U})^H\cdot\wbmi{U}=\mathbf{I}_{n_1},\quad \forall 1\leq i\leq n_3,\nonumber\\
& (\wbmi{V})^H\cdot\wbmi{V}=\mathbf{I}_{n_2},\quad \forall 1\leq i\leq n_3,\nonumber\\
& P_{\Omega'}(\wbmi{S}_t) = \wbmi{S}_t,~~\forall 1\leq i\leq
 n_3,1\leq t\leq T,\nonumber
\end{alignat}
where $\Omega'=\{(i,i):1\leq i\leq r \}$ is the set consisting of the indicators of all diagonal elements.

The alternative minimization algorithm proposed for solving Problem \eqref{fft_optim_problem} is described as follows.

\textbf{Update $\widehat{\bms}_t$}:
Problem \eqref{fft_optim_problem} with respect to $\wbmi{S}_t$ is:
\begin{alignat}{2} \label{subproblem:update_S}
    \mathop{\min}_{\{\wbmi{S}_t\}} &\sum_{t=p+1}^{T}\lVert\wbmi{S}_t- \sum_{i=1}^pa_i\wbmi{S}_{t-i}\rVert^2_F\nonumber\\
    &+\varphi\sum_{t=1}^{T}\lVert\wbmi{X}_t-\wbmi{U}\cdot\wbmi{S}_t\cdot(\wbmi{V})^H\rVert_F^2\nonumber\\
    \mbox{s.t.}\quad & P_{\Omega'}(\wbmi{S}_t) = \wbmi{S}_t,\quad \forall 1\leq t\leq T.
\end{alignat}

By computing the partial derivative of this cost function with respect to $\wbmi{S}_t$ and setting it equal to zero, we update $\wbmi{S}_t$ by
\begin{align}\label{update_S}
    \wbmi{S}_t = \begin{cases}
        \frac{1}{1+\varphi}P_{\Omega'}\left( \sum\limits_{i=1}^pa_i\wbmi{S}_{t-i}+\varphi (\wbmi{U})^H\wbmi{X}_t\wbmi{V}\right), \quad t>p,\\
        P_{\Omega'}\left((\wbmi{U})^H\wbmi{X}_t\wbmi{V}\right),\quad t\leq p.
    \end{cases}
\end{align}
As $i$ traverses from 1 to $n_3$, $\widehat{\bms}_t$ is completely updated.

\paragraph{Discussion 1: Relaxed-Diagonalization} We empirically explore the effect of relaxing the f-diagonalization in \eqref{subproblem:update_S}, {i.e.}, removing the constraint in Eq. \eqref{subproblem:update_S}.
One plausible explanation for the relaxed diagonal constraint is that even though we assumed in \eqref{X_tsvd} that the tensor time-series data has a shared factor subspace, in practice it tends to differ slightly, leading to the fact that part of the kernel tensor $\bms_t$ is not strictly diagonal.
Preserving small non-zero values on non-diagonal elements of $\bms_t$ may make the proposed model more flexible and richer in terms of the key information extracted from $\bmx_t$.
Under the Relaxed-Diagonalization setting, the update formula of $\wbmi{S}_t$ is
\begin{align}\label{update_S_relaxed}
    \wbmi{S}_t = \begin{cases}
        \frac{1}{1+\varphi}\left( \sum_{i=1}^pa_i\wbmi{S}_{t-i}+\varphi(\wbmi{U})^H\wbmi{X}_t\wbmi{V}\right),\quad t>p,\\
        (\wbmi{U})^H\wbmi{X}_t\wbmi{V},\quad t\leq p.
    \end{cases}
\end{align}

\textbf{Update $\widehat{\bmu},\widehat{\bmv}$}:
Problem \eqref{fft_optim_problem} with respect to $\wbmi{U}$ is:
\begin{alignat}{2}
    \mathop{\min}_{\wbmi{U}}&
    \sum_{t=1}^{T}\lVert\wbmi{X}_t-\wbmi{U}\cdot\wbmi{S}_t\cdot(\wbmi{V})^H\rVert_F^2\nonumber\\
    &\mbox{s.t.}\quad (\wbmi{U})^H\wbmi{U}=\mathbf{I}_{n_1}.
\end{alignat}
Then the update of $\wbmi{U}$ can be reduced to the following problem:
\begin{alignat}{2}\label{update_U_problem}
    \wbmi{U} &= \mathop{\arg\min}_{\wbmi{U}\in\mathbb{C}^{n_1\times r}} \sum_{t=1}^{T}\lVert\wbmi{X}_t-\wbmi{U}\cdot\wbmi{S}_t\cdot(\wbmi{V})^H\rVert_F^2\\
    &\overset{a}{=}\mathop{\arg\min}_{\wbmi{U}\in\mathbb{C}^{n_1\times r}}\sum_{t=1}^T\text{tr}\left(-\wbmi{U}\cdot\wbmi{S}_t\cdot\left(\wbmi{V}\right)^H\cdot\left(\wbmi{X}_t\right)^H \right)\nonumber\\
    &=\mathop{\arg\max}_{\wbmi{U}\in\mathbb{C}^{n_1\times r}}\text{tr}\left(\wbmi{U}\cdot\left(\sum_{t=1}^T\wbmi{X}_t\wbmi{V}\left(\wbmi{S}_t\right)^H \right)^H \right),\nonumber
\end{alignat}
where the equality $a$ holds due to $(\wbmi{U})^H\wbmi{U}=\textbf{I}_{n_1}$ and $\text{tr}(\cdot)$ denote the trace of matrix. Denote the SVD of $\sum_{t=1}^T\wbmi{X}_t\wbmi{V}\left(\wbmi{S}_t\right)^H$ as
\begin{equation}
    \sum_{t=1}^T\wbmi{X}_t\wbmi{V}\left(\wbmi{S}_t\right)^H = \bm{L_U^{i}}\text{diag}(\lambda_1,\ldots,\lambda_{n_1})\bm{{R_U^{i}}^H},
\end{equation}
 then the solution of problem \eqref{update_U_problem} can be expressed as
 \begin{equation}\label{update_U}
     \wbmi{U} = \bm{L_U^{i}}\bm{{R_U^{i}}^H}.
 \end{equation}
 Similarly, the update formula of $\wbmi{V}$ is shown as follows:
\begin{equation}\label{update_V}
    \wbmi{V} = \bm{L_V^{i}}\bm{{R_V^{i}}^H},
\end{equation}
where $\bm{L_V^{i}}$ and $\bm{R_V^{i}}$ are the left and right singular matrices of $\sum\limits_{t=1}^T\left(\wbmi{X}_t\right)^H\wbmi{U}\wbmi{S}_t$.

\textbf{Update $\{a_i \}$}:
The estimation of $\{a_i \}$ relies on the famous Yule-Walker equation \cite{yule1927vii}. We adopt a least squared modified Yule-Walker approach to update $\{a_i \}$ from the core tensors $\{\bms_t \}$.

\subsection{Step 2: Forecasting $\bmx_{T+1}$}
In this step, we first forecast the Fourier-transformed core tensor $\widehat{\bms}_{T+1}$ by the model learned in Step 1:
\begin{equation}\label{forecast S}
    \widehat{\bms}_{T+1} = \sum_{i=1}^pa_i\widehat{\bms}_{T+1-i}.
\end{equation}
After obtaining $\widehat{\bms}_{T+1}$, we apply IFFT to get $\bms_{T+1}=\mathcal{H}^{-1}(\widehat{\bms}_{T+1})$. Finally, we reconstruct $\bmx_{T+1}$ by t-SVD model with optimized orthogonal tensors $\bmu$ and $\bmv$:
\begin{equation}\label{forecast X}
    \bmx_{T+1} = \bmu\ast\bms_{T+1}\ast\bmv^H.
\end{equation}

We summarize the whole \textbf{LOTAP} process in \textbf{Algorithm 1}. {Note that while truncated t-SVD is introduced for modeling purposes in the optimization formulation, no explicit t-SVD computations are involved in the actual iterative procedure. This is also the key reason why our algorithm enjoys an efficient per-iteration complexity.}

\begin{algorithm}[ht]
\caption{LOTAP for time-series forecasting problem}\label{algorithm:LOTAP}
\begin{algorithmic}[1]
\renewcommand{\algorithmicrequire}{\textbf{Input:}}
\renewcommand{\algorithmicensure}{\textbf{Output:}}
\REQUIRE A time-series data $\bmx\in\mathbb{C}^{n_1\times n_2\times n_3\times T}$ and core size $r$;
\ENSURE $\bmx_{T+1},\bmu,\bmv$.

\STATE \underline{\textbf{Step 1: Tensor AR with truncated t-SVD}}
\STATE Initialize $\bmu$ and $\bmv$ randomly.
\STATE Initialize $\bms_t$ by $\bms_t=\bmu^H\ast\bmx_t\ast\bmv$ for $t=1,\dots,T.$
\STATE $\widehat{\bms}_t\leftarrow \mathcal{H}(\bms_t),\quad  t=1,\dots,T.$
\STATE $\widehat{\bmu}\leftarrow \mathcal{H}(\bmu),\widehat{\bmv}\leftarrow \mathcal{H}(\bmv).$
\WHILE{ not convergence }
    \STATE  Estimate coefficients $\{a_i \}_{i=1}^p$ of \textbf{AR} via Yule-Walker equations based on $\{\bms_t \}_{t=1}^T.$
    \STATE Update $\wbmi{S}_t$ by \eqref{update_S} for $t=1,\dots,T$ and $i=1,\dots,n_3$.
    \STATE If applying relaxed-diagonalization, update $\wbmi{S}_t$ by \eqref{update_S_relaxed} for $t=1,\dots,T$ and $i=1,\dots,n_3$.
    \STATE Update $\wbmi{U}$ by \eqref{update_U} for $i=1,\dots,n_3$.
    \STATE Update $\wbmi{V}$ by \eqref{update_V} for $i=1,\dots,n_3$.
\ENDWHILE
\STATE \underline{\textbf{Step 2: Forecasting $\bmx_{T+1}$}}
\STATE Estimate the Fourtier-transformed core tensor $\widehat{\bms}_{T+1}$ by \eqref{forecast S}.
\STATE $\bmu\leftarrow \mathcal{H}^{-1}(\widehat{\bmu}),\bmv\leftarrow \mathcal{H}^{-1}(\widehat{\bmv}).$
\STATE Compute $\bmx_{T+1}$ by \eqref{forecast X}.
\RETURN $\bmx_{T+1},\bmu,\bmv$.
\end{algorithmic}
\end{algorithm}

{\paragraph{Discussion 2: Extension to Time-Series Imputation}  
The core idea of LOTAP can, in principle, be extended to handle time-series imputation by leveraging observations both before and after the missing interval. This would require introducing a masking operator in the data-fidelity term, under which the updates for $\bmu, \bmv, \bms$ would remain largely similar with only minor modifications. The main challenge arises in updating the autoregressive parameters $\alpha_i, \beta_i$: while our current method relies on Yule–Walker equations that assume complete data, missing entries would lead to biased covariance estimates and thus invalidate the closed-form updates. Addressing this issue may involve EM-based maximum likelihood or related estimation techniques. We leave a rigorous development of this imputation extension as an interesting direction for future work.}

{\paragraph{Discussion 3: Direction-dependence of t-SVD.}  
It is worth noting that t-SVD is direction-dependent since the Fourier transform is applied along the third mode. In practice, this is often not a drawback but rather an advantage, as the third dimension is typically chosen to capture periodic structures. For example, in our experiments on the USHCN dataset, the first two dimensions represent spatial locations (latitude and longitude), while the third dimension corresponds to one week (five days), which naturally exhibits periodicity. Moreover, the concept of multi-tubal rank has been proposed in the literature to address this limitation by performing decompositions along different modes \cite{tsvd-tucker-proof2020zhang,qin2025low}, and our LOTAP framework can, in principle, be extended in that direction. We leave such extensions as an interesting direction for future work.}

\section{Complexity Analysis}\label{Subsec: Complexity}

Given that prior time series forecasting algorithms centered around tensor structures, such as MOAR \cite{moar2018jing}, MCAR \cite{moar2018jing}, and BHT-ARIMA \cite{bhtarima2020shi}, are grounded in tensor Tucker decomposition, this section is first dedicated to dissecting the distinctions between this category of algorithms and LOTAP based on the truncated t-SVD. The another purpose is to establish complexity analysis for LOTAP.

The class of algorithms based on t-SVD presents a wider range of applicability when contrasted with the algorithms grounded in Tucker decomposition. Both categories of algorithms necessitate the underlying time series data to adhere to a low-rank structure as per their respective decomposition methodologies. A proof presented in \cite{tsvd-tucker-proof2020zhang} establishes that: $$\text{rank}_t(\bmx)\leq \min\{\text{rank}(\boldsymbol{X}_{(1)}), \text{rank}(\boldsymbol{X}_{(2)})\},$$ where $\boldsymbol{X}_{(i)}$ is the mode-$i$ matricization of $\bmx$. This equation illustrates that a third-order tensor characterized by low-rank properties, as stipulated by the Tucker decomposition framework, inherently retains its low-rank nature within the t-SVD paradigm, and vice versa. As a result, the applicability of our LOTAP algorithm extends across the domain of previously employed Tucker decomposition algorithms. Notably, our algorithm also tackles specific instances of time series data that posed challenges for the original algorithms' efficacy. A comprehensive elaboration of this proposition will be provided in the ensuing experimental section.

The LOTAP algorithm demonstrates superior performance in terms of computational speed compared to the algorithm grounded in Tucker decomposition. This discrepancy stems from the inherent dissimilarity in the computational efficiency of the two distinct decomposition methodologies. The Tucker decomposition necessitates SVD of the complete tensor-reconstructed matrix during computation. Conversely, t-SVD computation entails separate SVD for each frontal slice of the tensor following FFT. Notably, practical implementation reveals that the Tucker decomposition requires approximately ten times the duration consumed by the t-SVD. We will empirically validate this assertion in the subsequent numerical experiment section.

There exists no pronounced advantage or drawback concerning prediction accuracy. Theoretically, the kernel tensor obtained by t-SVD assumes an f-diagonal structure, while the kernel tensor derived by Tucker decomposition contains a higher information density than its t-SVD counterpart of equivalent dimensions. This phenomenon suggests that algorithms modeled after Tucker decomposition, which leverage this richer information, could potentially exhibit superior predictive performance compared to the LOTAP algorithm. However, empirical evidence derived from numerical experiments contradicts this expectation, revealing negligible differences in predictive accuracy between the two methodologies.

We now conduct a thorough examination of the computational complexity associated with Algorithm \ref{algorithm:LOTAP}. We present the detailed computational cost of each step within a single iteration, which can be found in Table \ref{table: complexity of LOTAP}. Additionally, we offer a comparative analysis with other algorithms \cite{bhtarima2020shi,moar2018jing} in Table \ref{table: complexity comparison}. This analysis will shed light on the efficiency and performance of our proposed algorithm LOTAP in relation to existing methods.

At each iteration, the computational cost of the Yule-Walker method amounts to $O(p^3+pTrn_3)$. Subsequently, the complexity associated with updating $\{\bms_t\}_{t=1}^T$ using Eq. \eqref{update_S} is $O(pTn_3r+2Tn_1n_2)$. Moreover, we employ SVD to update $\bmu$ and $\bmv$ following Eqs. \eqref{update_U} and \eqref{update_V}, respectively. The computational cost for these updates is $O((n_1+n_2)n_3r^2+2Tn_1n_2n_3)$. Consequently, the overall computational complexity for each iteration is given by $O(Tn_1n_2n_3+(n_1+n_2)n_3r^2)$.

Upon comparison, we find that the computational costs for MOAR and MCAR algorithms in \cite{moar2018jing} are respectively $O(Trn_1n_2n_3+n_1n_2n_3r)$ and $O(Trn_1n_2n_3+(n_1^3+n_2^3+n_3^3)r)$ during each iteration, where the Tucker rank of $\bmx_t$ is represented as $(r,r,r)$. Furthermore, BHT-ARIMA (\cite{bhtarima2020shi}), which builds upon MCAR via utilizing the MDT strategy and ARIMA model, incurs a cost of $O(Tn_1n_2n_3\tau r)$ at each iteration, where $\tau$ denotes the MDT length.

Evidently, our proposed method demonstrates superior efficiency compared to state-of-the-art algorithms. This advantage arises from the approach of breaking down the SVD of a large matrix into individual SVDs of $n_3$ frontal slice matrices. Consequently, our algorithm achieves more efficient updates for the projection matrix and then outperforms Tucker-decomposition-based algorithms, which necessitate the computation of the SVD of a larger matrix. This key difference in computational strategy contributes significantly to the improved efficiency of our approach.

\begin{table}[h!]
    \centering
    \caption{The per-iteration complexity of LOTAP}
    \begin{tabular}{c c}
    \hline
    Step & Computational Complexity\\ \hline
    update $\{a_i\}_{i=1}^p$ & $O(p^3+pTrn_3)$\\
    update $\bmu$ & $O(Tn_1n_2n_3+n_1n_3r^2)$\\
    update $\bmv$ & $O(Tn_1n_2n_3+n_2n_3r^2)$\\
    update $\{\bms_t\}_{t=1}^T$ & $O(pTn_3r+Tn_1n_2)$\\ \hline
    Total & $O(Tn_1n_2n_3+(n_1+n_2)n_3r^2)$\\ \hline
    \end{tabular}
    \label{table: complexity of LOTAP}
\end{table}

\begin{table}[h!]
    \centering
    \caption{Complexity comparison of several algorithms}
    \begin{tabular}{c c}
    \hline
    Algorithm & Computational Complexity\\ \hline
    LOTAP & $O(Tn_1n_2n_3+(n_1+n_2)n_3r^2)$\\
    MOAR \cite{moar2018jing} & $O(Trn_1n_2n_3+n_1n_2n_3r)$\\
    MCAR \cite{moar2018jing} & $O(Trn_1n_2n_3+(n_1^3+n_2^3+n_3^3)r)$\\
    BHT-ARIMA \cite{bhtarima2020shi} & $O(Tn_1n_2n_3\tau r)$\\ \hline
    \end{tabular}
    \label{table: complexity comparison}
\end{table}

\section{Numerical Experiments}\label{Section: Numerical}
In this section, we present the validation of our proposed algorithm, LOTAP, for modeling higher-order time series in the context of time series forecasting. We conduct experiments using one synthetic dataset and three real-world datasets, with sizes ranging from approximately 100 to approximately 1000 data points. We provide a concise overview of the datasets and experimental settings before presenting the forecasting results.
All numerical experiments are performed using MATLAB R2022b on a Windows PC equipped with a 14-core Intel(R) Core(TM) 2.30GHz CPU and 16GB RAM.

\subsection{Datasets and Experimental Settings}
For our time series forecasting  experiments, we utilize the following datasets, which encompass one synthetic dataset (SYN) and three real-world datasets (USHCN, NASDAQ100, CCDS). These datasets serve as the basis for evaluating the performance and effectiveness of our LOTAP algorithm\footnote{Available at https://github.com/whn18/LOTAP}.

\textit{SYN}: The synthetic (SYN) dataset is a low-rank, third-order tensor time series that we generated using the following method. We first generate the core f-diagonal tensor series $\{\bms_t \}_{t=1}^{1000}$ with AR(3) model. Then we generate $10$ random matrices $\boldsymbol{A}_1,\dots,\boldsymbol{A}_{10}\in\mathbb{C}^{100\times 4}$ with \textit{i.i.d} $\mathcal{N}(0,1)$ entries. Let $\widehat{\boldsymbol{U}}^{(1)},\dots,\widehat{\boldsymbol{U}}^{(10)}$ be orthonormal bases for their respectively column spaces and then generate $\widehat{\bmu}\in\mathbb{C}^{100\times 4\times 10}$ whose $k$-st frontal slice is $\widehat{\boldsymbol{U}}^{(k)}$. By using Matlab command $\texttt{ifft}(\widehat{\bmu},[],3)$, we obtain column-orthogonal tensor $\bmu\in\mathbb{C}^{100\times 4\times 10}$. $\bmv\in\mathbb{C}^{100\times 4\times 10}$ is generated by the same process as $\bmu$. The generated noisy tensor time series are formulated as
\begin{equation}
    \bmx_t = \bmu*\bms_t*\bmv^H+\rho \lVert\bms_t\rVert_F\bme_t\in\mathbb{C}^{100\times100\times10},
\end{equation}
where $\rho=0.01$ is the noise parameter and $\bme_t$ is the noise tensor with \textit{i.i.d.} $\mathcal{N}(0,1)$ entries. It is notable that we have $\lVert S_t\rVert_F = \lVert \bmu\ast\bms_t\ast\bmv^H\rVert_F$. The successful outcome of this process results in the generation of a synthetic time series dataset with dimensions of $100\times 100\times 10$, spanning a total of 1000 time points.

\textit{USHCN}\footnote{https://www.ncei.noaa.gov/pub/data/ushcn/v2.5/}: The U.S. historical climate network (USHCN) dataset records climate data for each month at various locations across the United States. It includes four climate statistics: average monthly maximum temperature, average monthly minimum temperature, average monthly air temperature, and total monthly precipitation. For our experiment, we focus on quarterly average observations from 120 weather stations over a period of 75 years, spanning from 1940 to 2014. The data is organized into a $120\times 4\times 4\times 75$ time series.
\begin{enumerate}
    \item[1)] The first dimension (120) represents the 120 weather stations.
    \item[2)] The second dimension (4) denotes the four climate statistics (average monthly maximum temperature, average monthly minimum temperature, average monthly air temperature, and total monthly precipitation).
    \item[3)] The third dimension (4) is used to represent the four quarters of each year (75 time points in total).
    \item[4)] The fourth dimension (75) represents the 75 years of data.
\end{enumerate}
This arrangement allows us to effectively analyze and forecast climate patterns across the United States based on the quarterly average observations collected from the \textit{USHCN} dataset.

\textit{NASDAQ100}\footnote{https://github.com/alireza-jafari/GCNET-Dataset}: The Nasdaq 100 Index is a basket of the 100 largest, most actively traded companies listed on the Nasdaqstock exchange. The NASDAQ100 dataset \cite{nasdaq2022} includes five key indicators, namely opening, closing, adjusted-closing, high,and low prices and volume. In this paper, we use five consecutive years of daily data for 100 companies from 2016 onwards, reconstructed by week, to obtain a $100\times5\times5\times250$ times series.
\begin{enumerate}
    \item[1)] The first dimension (100) represents the 100 companies.
    \item[2)] The second dimension (5) denotes the five key index (opening, closing, adjusted-closing, high,and low prices and volume).
    \item[3)] The third dimension (5) is used to represent the five days of the week when stocks are open.
    \item[4)] The fourth dimension (250) represents the 250 weeks of data.
\end{enumerate}

\textit{CCDS}\footnote{https://melady.usc.edu/data/}: The Comprehensive Climate Dataset (CCDS) is a collection of climate records of North America from \cite{CCDS2009}. In this paper, the data are 2.5-by-2.5 grids from $30^{\circ}$N,  $92^{\circ}$E to $45^{\circ}$N,  $107^{\circ}$E recording 17 variables such as carbon dioxide monthly from 1990 to 2001, yielding 144 epochs.

In our experiments, we ensure a rigorous evaluation of our LOTAP algorithm by dividing all datasets into training and test sets. For the SYN dataset, the first 80 time points are used to train the models, and the remaining time points constitute the test set. For the USHCN dataset, we use the first 40 time points to train the models and forecast the following 35 time points. For the NASDAQ100 dataset, we train the first 200 epochs and forecast the next 50 epochs. For the CCDS dataset, we use the first 120 time points to forecast the remaining 24 time points.

In our comparative study, we evaluated five competing methods for time series forecasting: 1) the classical \textbf{ARIMA}; 2) the popular industrial forecasting method: Amazon-\textbf{DeepAR} \cite{salinas2020deepar}; 3) the three tucker-decomposition-based methods: \textbf{MOAR}, \textbf{MCAR} \cite{moar2018jing} and \textbf{BHT-ARIMA} \cite{bhtarima2020shi}. To ensure a fair comparison, we performed a large number of experiments on the synthetic dataset and multiple repetitions on the real-world datasets. For each algorithm we performed 1000 Monte-Carlo experiments on the synthetic dataset and 10 repetitions on the real-world datasets.

In the development phase, we use a grid search to find the optimal parameters for the LOTAP algorithm. The maximum number of iterations is set to 30, serving as the stopping criterion for the algorithm to prevent excessive computation.

To evaluate the forecasting accuracy of the LOTAP algorithm and compare it with other methods, we utilize the widely used Mean Squared Percentage Error (MSPE) metric. The MSPE at time point $t$ is computed by
\begin{equation}
    \text{MSPE} = \mathbb{E}\left(\frac{\lVert \bmx_t-\bm{\mathcal{P}}_t \rVert_F}{\lVert\bmx_t \rVert_F}\right),
\end{equation}
where $\bm{\mathcal{P}}_t$ is the forecasted tensor of $\bmx$ at time point $t$ and $\mathbb{E}(\cdot)$ is the expectation with respect to $t$.
In our experiments, after successive predictions over the test sets, we calculate the relative error at each prediction point and average then to obtain the value of MSPE.
The MSPE metric quantifies the relative forecasting error with respect to the true tensor values at each time point. Smaller MSPE values indicate better forecasting performance, as they imply that the forecasted tensor is closer to the ground truth data. This metric allows us to assess and compare the accuracy of different methods in forecasting time series data effectively.

\subsection{Low-rank Verification}
In this subsection, we will verify that the four datasets (SYN, HSHCN, NASDAQ100, CCDS) do actually have a low-tubal-rank structure under the Fourier transform. Furthermore, we will show that these data exhibit a stronger low-dimensional structure under t-SVD compared to Tucker decomposition.

\begin{figure}[!t]
    \centering
    \subfigure[]{
        \begin{minipage}[t]{0.49\linewidth}
            \centering
            \includegraphics[width=0.9\linewidth]{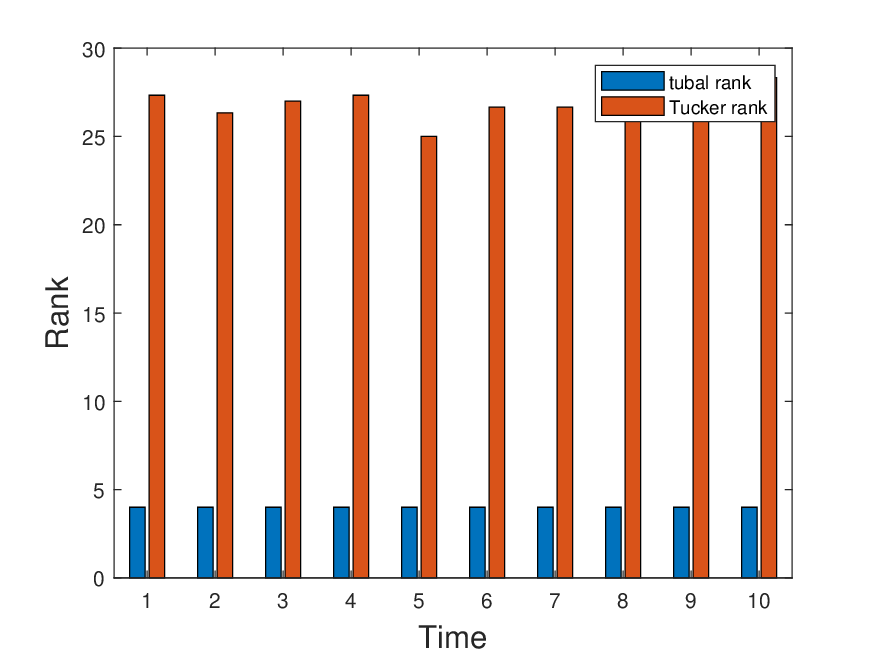}
        \end{minipage}
    }%
    \subfigure[]{
        \begin{minipage}[t]{0.49\linewidth}
            \centering
            \includegraphics[width=0.9\linewidth]{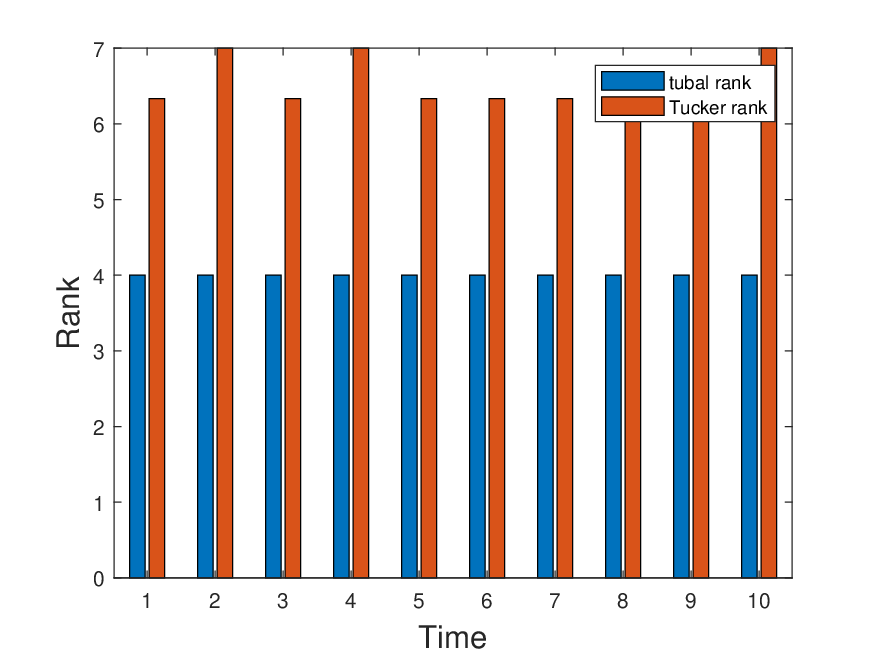}
        \end{minipage}
    }%

    %´Ë´¦µÄ¿ÕÐÐºÜÖØÒª£¬ÏëÈÃÍ¼Æ¬ÔÚÊ²Ã´µØ·½»»ÐÐ¾ÍÔÚ´úÂë¶ÔÓ¦Î»ÖÃ¿ÕÐÐ
    \subfigure[]{
        \begin{minipage}[t]{0.49\linewidth}
            \centering
            \includegraphics[width=0.9\linewidth]{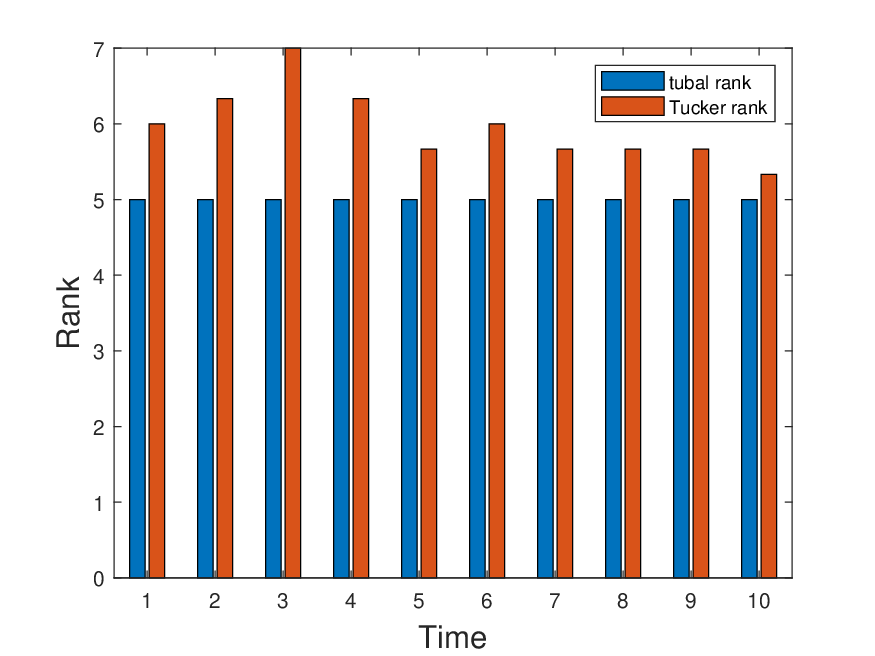}
        \end{minipage}
    }%
    \subfigure[]{
        \begin{minipage}[t]{0.49\linewidth}
            \centering
            \includegraphics[width=0.9\linewidth]{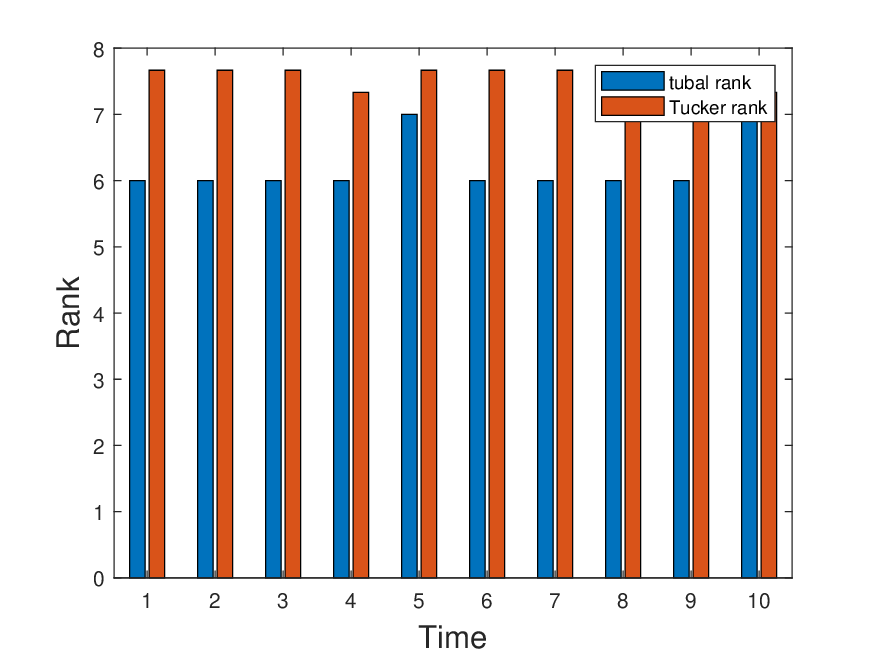}
        \end{minipage}
    }%
    \centering
    \caption{\small{Comparison of tubal rank with average Tucker rank on the a) SYN dataset, b) USHCN dataset, c) NASDAQ100 dataset, d) CCDS dataset.}}
    \label{Fig: rank}
\end{figure}

To demonstrate this, we select $10$ successive time points in these datasets and examine both their tubal rank and tucker rank. Due to the presence of noise in the data, we approximate those singular values that are less than $10^{-2}$ times the main singular value as 0 when calculating both the tubal rank and Tucker rank. As well, for comparison purposes, we consider the so-called average Tucker rank, which is the average of the ranks of the three modes, as follows
\begin{equation}
    \text{rank}_{\text{avg}}(\bmx) = \frac{\text{rank}(\boldsymbol{X}_{(1)})+\text{rank}(\boldsymbol{X}_{(2)})+\text{rank}(\boldsymbol{X}_{(3)})}{3},
\end{equation}
where $\boldsymbol{X}_{(i)}$ is the mode-$i$ matricization of $\bmx$. The experimental results are shown in Figure \ref{Fig: rank}. For SYN dataset, it is apparent that our synthetic data does not exhibit low-dimensional structure under Tucker decomposition, whereas it is low-rank under t-SVD.
As for three datasets, they all exhibit a lower rank structure under t-SVD compared to Tucker decomposition. Our experimental results validate our theoretical analysis in Section \ref{Subsec: Complexity}: third-order tensors with low-rank properties, as specified by the Tucker decomposition framework, inherently retain their low-rank properties in the t-SVD paradigm.

\subsection{Verification of Subspace Stability Across Time} 
{In this subsection, we validate the assumption stated at the beginning of Section~\ref{Section: Alg-des}, namely that the subspace factors of adjacent tensor data along the time scale vary only slightly. To this end, we evaluate
\[
	\mathrm{res}_t := \frac{1}{\sqrt{2\min\{n_1,n_2\}}}\min_{\Theta}\lVert\mathcal{U}_1 - \mathcal{U}_t \ast \Theta\rVert_F,
\]
where $\Theta$ is a column-orthogonal tensor, and $\mathcal{U}_1$ and $\mathcal{U}_t$ are the left factor tensors obtained from the truncated t-SVD of $\mathcal{X}_1$ and $\mathcal{X}_t$, respectively. The normalization factor $\frac{1}{\sqrt{2\min\{n_1,n_2\}}}$ ensures that $\mathrm{res}_t \in [0,1]$. This quantity measures the angular discrepancy between the subspaces spanned by $\mathcal{U}_1$ and $\mathcal{U}_t$. 
We omitted here the experiment of $\mathcal{V}$ since its behavior is analogous to that of $\mathcal{U}$.
In the experiments, we compute $\mathrm{res}_t$ using the first 50 time points from each of the four datasets. During the truncated t-SVD of $\mathcal{X}_1$ and $\mathcal{X}_t$, singular values smaller than $0.01$ of the largest singular value are discarded, and the corresponding columns in $\mathcal{U}$ are set to zero vectors.
}

\begin{figure}[ht]
    \centering
    \includegraphics[width=0.75\linewidth]{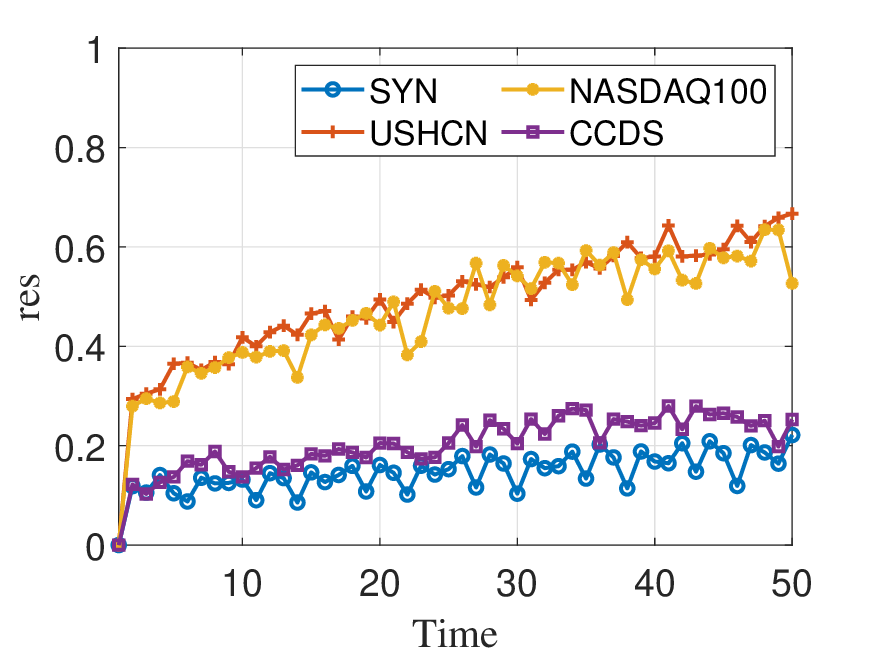}
    \caption{{Evolution of $\mathrm{res}_t$ over time for the four datasets, illustrating the degree of variation in subspace factors across adjacent time points.}}
    \label{fig:subspace_res}
\end{figure}

{As shown in Figure~\ref{fig:subspace_res}, $\mathrm{res}_t$ remains close to zero when $t$ is small and gradually increases as $t$ grows, which numerically confirms our assumption that subspace factors vary slowly across adjacent time points.}

\subsection{Performance Evaluation}
In this subsection, we evaluate our proposed algorithm in the following five aspects to show the out-performance of LOTAP.

\paragraph{1) Convergence criterion and maximum iteration number}
We investigate the convergence of the LOTAP algorithm by examining the relative error of the projection column-orthogonal tensors $\bmu$ and $\bmv$. The relative error can be computed as follows:
$$\frac{\lVert\bmu^{k+1}-\bmu^k\rVert_F^2 + \lVert\bmv^{k+1}-\bmv^k \rVert_F^2}{\lVert\bmu^{k+1} \rVert_F^2+\lVert\bmv^{k+1} \rVert_F^2},$$
where $\bmu^{k+1}$ and $\bmv^{k+1}$ represent the tensors obtained after the $(k+1)$-th iteration, and $\bmu^k$ and $\bmv^k$ represent the tensors obtained after the $k$-th iteration.

\begin{figure}[htbp]
    \centering
    \subfigure[]{
        \begin{minipage}[t]{0.49\linewidth}
            \centering
            \includegraphics[width=0.9\linewidth]{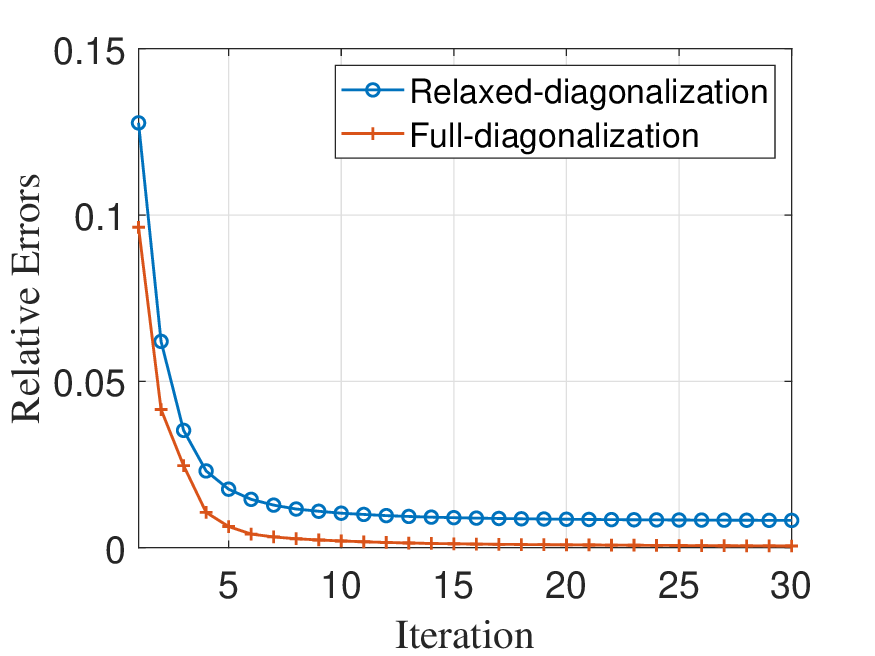}
        \end{minipage}
    }%
    \subfigure[]{
        \begin{minipage}[t]{0.49\linewidth}
            \centering
            \includegraphics[width=0.9\linewidth]{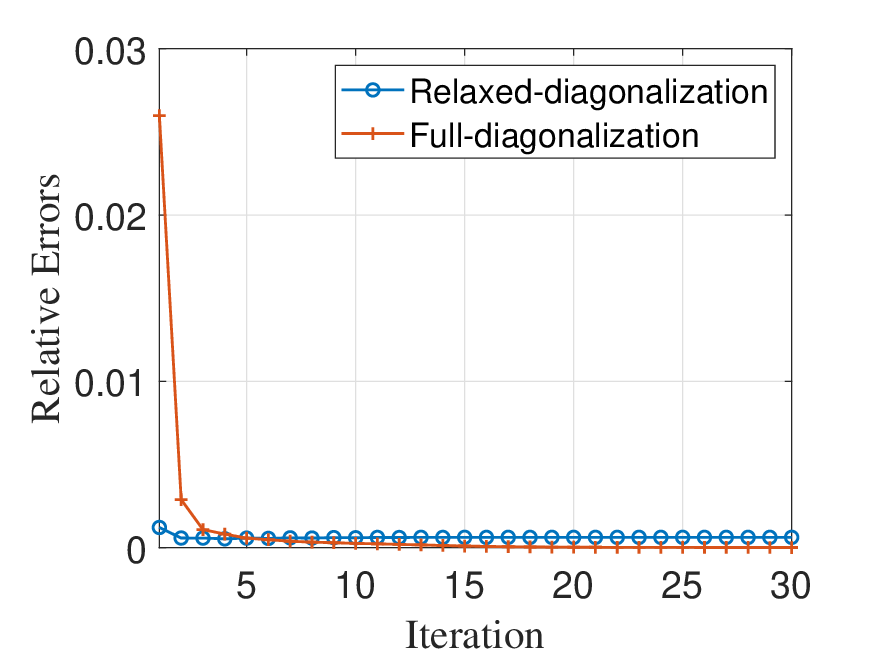}
        \end{minipage}
    }%

      \subfigure[]{
        \begin{minipage}[t]{0.49\linewidth}
            \centering
            \includegraphics[width=0.9\linewidth]{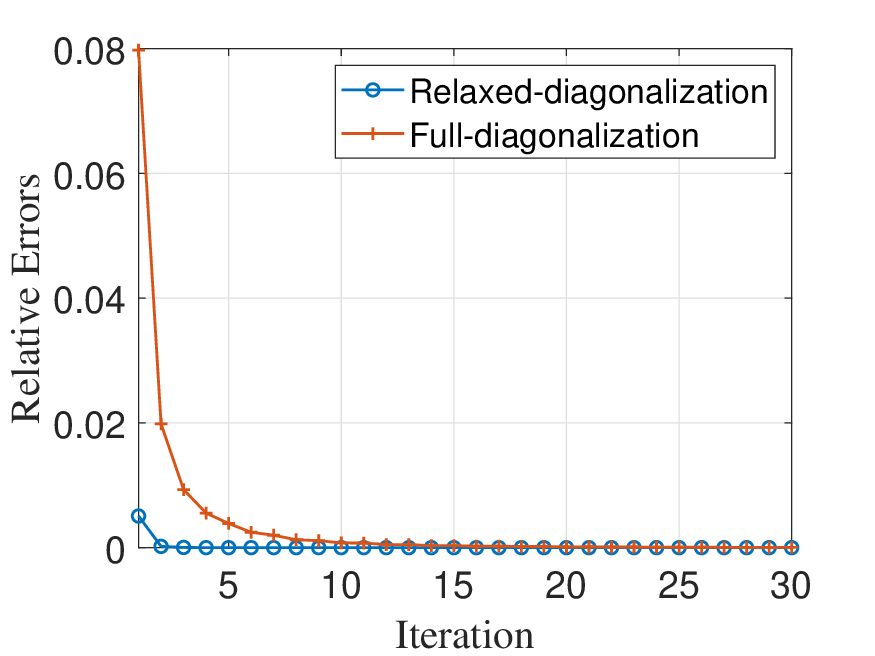}
        \end{minipage}
    }%
    \subfigure[]{
        \begin{minipage}[t]{0.49\linewidth}
            \centering
            \includegraphics[width=0.9\linewidth]{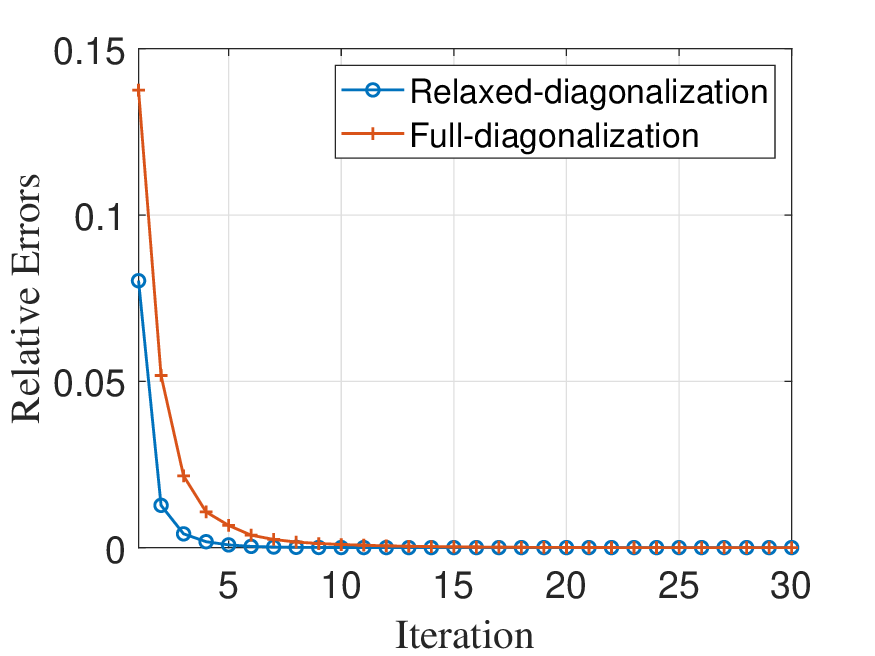}
        \end{minipage}
    }%
    \centering
    \caption{\small{Convergence curves of LOTAP on the a) SYN dataset, b) USHCN dataset, c) NASDAQ100 dataset, d) CCDS dataset.}}
    \label{Fig: Convergence}
\end{figure}

Fig. \ref{Fig: Convergence} displays the variation of the relative error of $\bmu$ and $\bmv$ during the iterations. As shown in the figure, the relative errors obtained by the LOTAP algorithm with the full-diagonalization and the relaxed-diagonalization  versions both decrease rapidly and converge in no more than 10 iterations. The iteration criterion plays a vital role in ensuring the convergence of the objective function. In this paper, we adopt a stopping criterion for our proposed method, where the relative change between two consecutive iterations must be below $1e-3$, and the maximum number of iterations is set to $10$. This criterion ensures the algorithm's convergence while effectively managing computational resources.

\paragraph{2) Sensitivity analysis of parameters}
We conducted a sensitivity analysis to investigate the impact of different parameter choices in the LOTAP algorithm. Note that LOTAP choose used in the following experiments is the Relaxed-diagonalization version. There are two crucial parameters in our LOTAP method: the core size $r$, the AR order $p$ and the regularization parameter $\varphi$. In our experiments, the parameters $p$ and $\varphi$ are selected via grid search in a heuristic manner. As for the core size $r$, we choose it naturally through the results of the t-SVD decomposition.

\begin{figure}[htbp]
    \centering
    \subfigure[]{
        \begin{minipage}[t]{0.49\linewidth}
            \centering
            \includegraphics[width=0.9\linewidth]{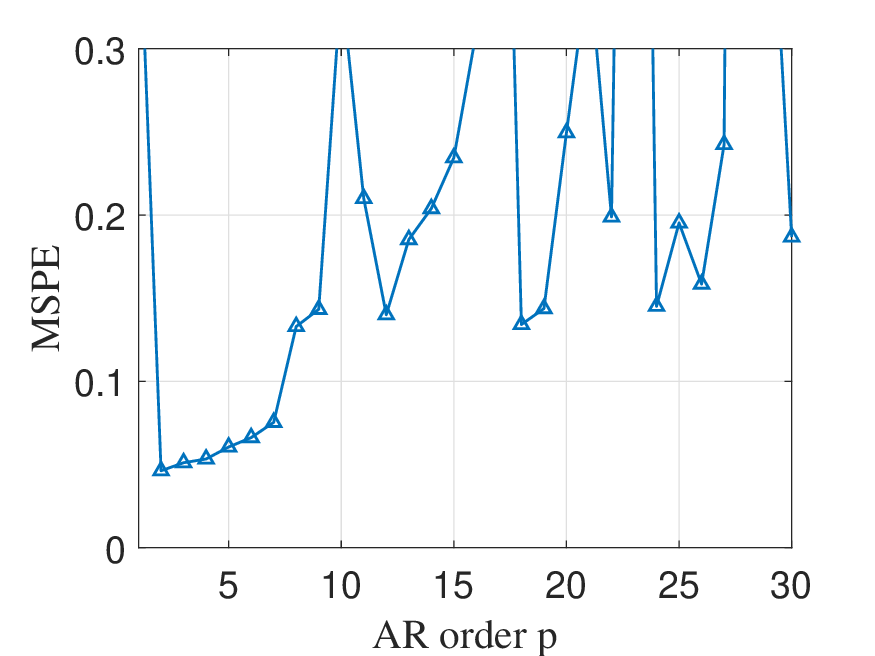}
        \end{minipage}
    }%
    \subfigure[]{
        \begin{minipage}[t]{0.49\linewidth}
            \centering
            \includegraphics[width=0.9\linewidth]{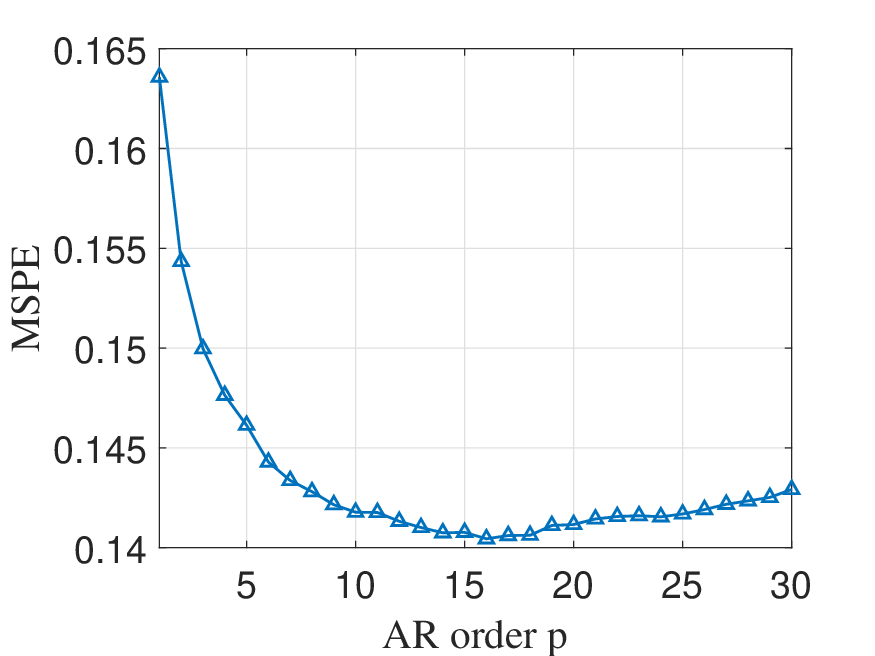}
        \end{minipage}
    }%

    \subfigure[]{
        \begin{minipage}[t]{0.49\linewidth}
            \centering
            \includegraphics[width=0.9\linewidth]{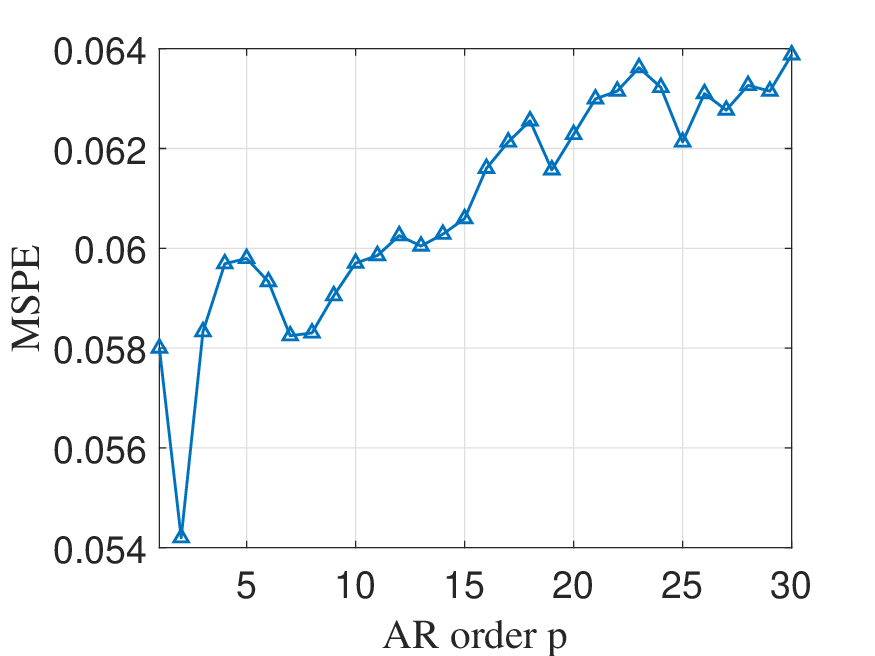}
        \end{minipage}
    }%
    \subfigure[]{
        \begin{minipage}[t]{0.49\linewidth}
            \centering
            \includegraphics[width=0.9\linewidth]{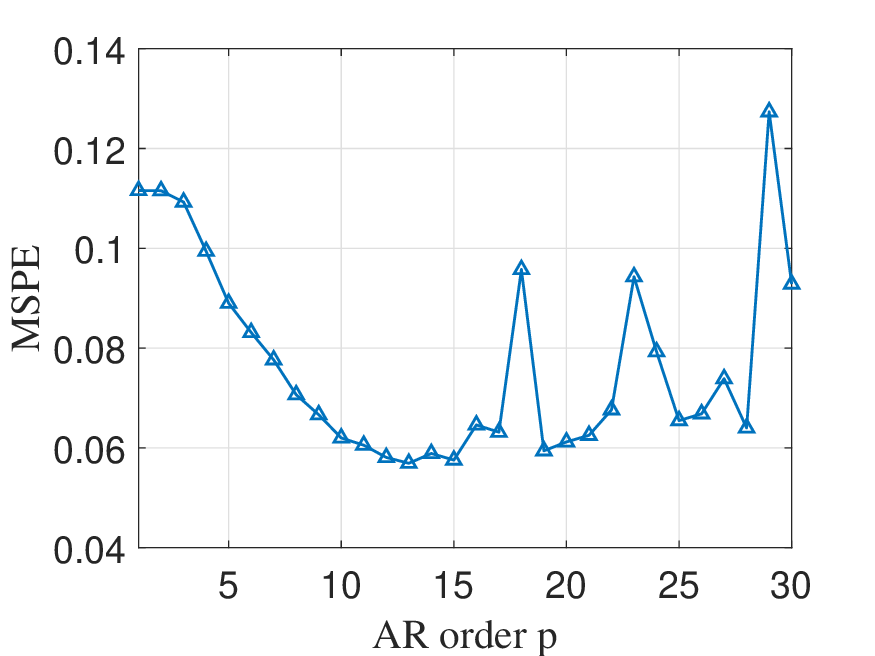}
        \end{minipage}
    }%
    \centering
    \caption{\small{Effect of the parameter $p$ for the LOTAP algorithm on the a) SYN dataset, b) USHCN dataset, c) NASDAQ100 dataset, d) CCDS dataset.}}
    \label{Fig: sensitivity-order}
\end{figure}

Fig. \ref{Fig: sensitivity-order} illustrates the average MSPE values obtained by the LOTAP algorithm for varying values of the AR order $p$. We found that the optimal AR orders for the four datasets are $\{2,16,2,13\}$. Notably, in the case of the synthetic dataset (see Fig. \ref{Fig: sensitivity-order} (a)), the MSPE values appear to deteriorate significantly when an inappropriate order $p$ is selected. In contrast, for the three real-world datasets (see Fig. \ref{Fig: sensitivity-order} (b)-(d)), the MSPE values do not exhibit such drastic degradation even with suboptimal AR orders. This observation highlights the robust accuracy of the LOTAP algorithm when faced with the selection of the AR order $p$ in real-world scenarios.

\begin{figure}[h!]
    \centering
    \subfigure[]{
        \begin{minipage}[t]{0.49\linewidth}
            \centering
            \includegraphics[width=0.9\linewidth]{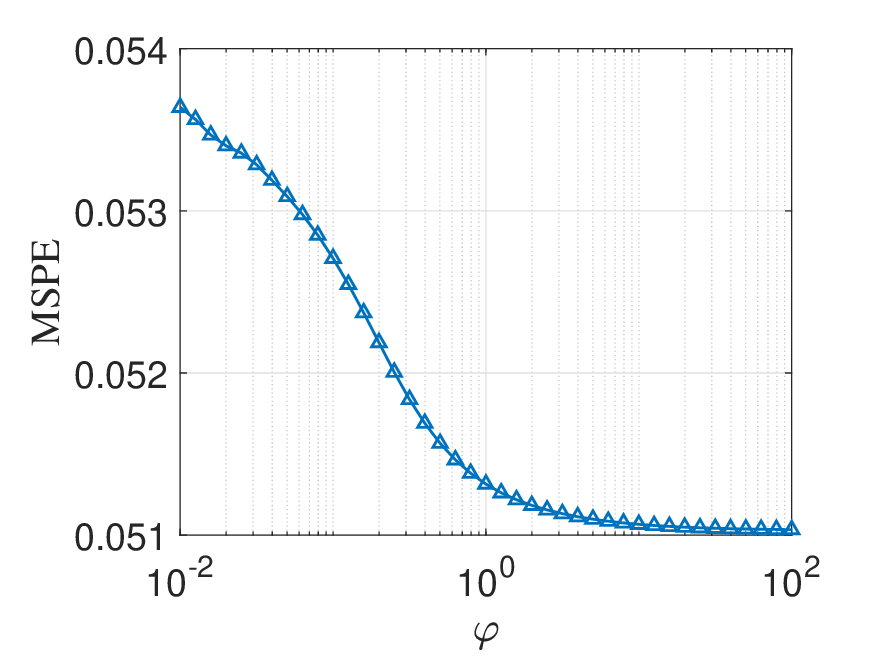}
        \end{minipage}
    }%
    \subfigure[]{
        \begin{minipage}[t]{0.49\linewidth}
            \centering
            \includegraphics[width=0.9\linewidth]{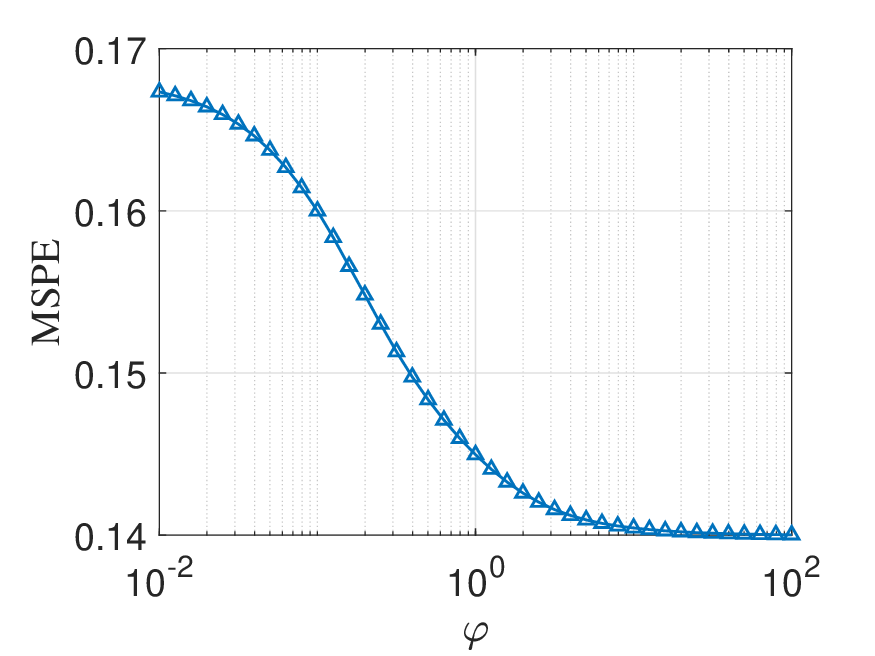}
        \end{minipage}
    }%

    \subfigure[]{
        \begin{minipage}[t]{0.49\linewidth}
            \centering
            \includegraphics[width=0.9\linewidth]{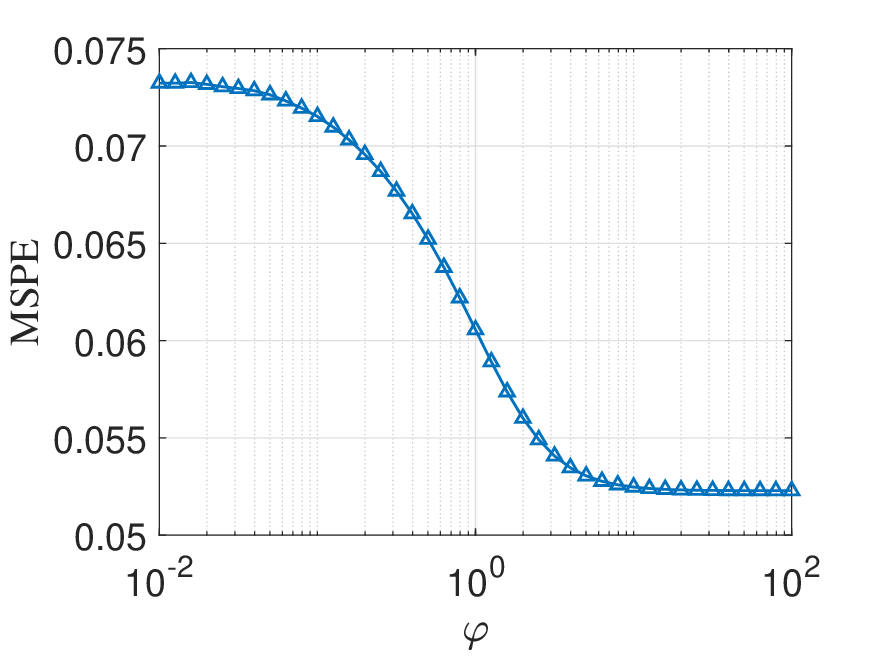}
        \end{minipage}
    }%
    \subfigure[]{
        \begin{minipage}[t]{0.49\linewidth}
            \centering
            \includegraphics[width=0.9\linewidth]{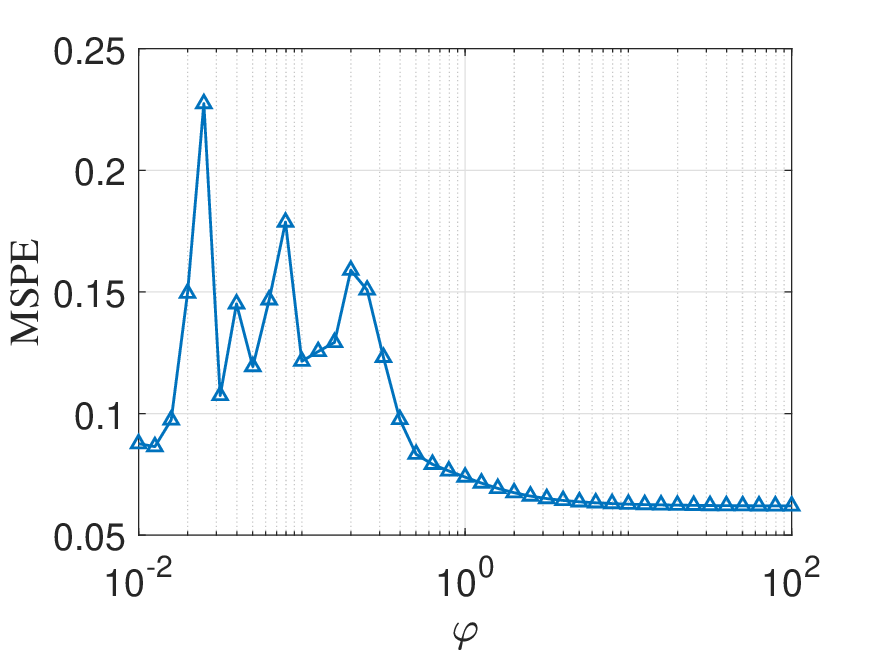}
        \end{minipage}
    }%
    \centering
    \caption{\small{Effect of the parameter $\varphi$ for the LOTAP algorithm on the a) SYN dataset, b) USHCN dataset, c) NASDAQ100 dataset, d) CCDS dataset. When $\varphi$ is large enough, the accuracy of LOTAP is not very sensitive to the choice of parameter $\varphi$.}}
    \label{Fig: sensitivity-phi}
\end{figure}

Fig. \ref{Fig: sensitivity-phi} illustrates the performance of the LOTAP algorithm across various values of the regularization parameter $\varphi$, spanning from $0.01$ to $100$. Upon observation, it becomes apparent that when smaller $\varphi$ values are utilized, the forecasting performance is suboptimal for all four datasets. This outcome can be attributed to the fact that smaller $\varphi$ values lead to the regularization term in Problem \eqref{fft_optim_problem} approaching zero. Consequently, the t-SVD becomes inaccurate, impeding the extraction of critical information from the original time tensor series.

Conversely, as $\varphi$ increases and reaches a sufficiently large value, the forecasting accuracy and robustness increase and show a tendency to converge. The regularization term plays a significant role in maintaining the balance between accurately modeling the data and avoiding overfitting. Upon the judicious selection of $\varphi$, the algorithm attains an equilibrium that fosters enhanced predictive performance while mitigating the occurrence of undue iterations within the model.

\paragraph{3) Comparison of forecasting accuracy}
We present the forecasting performance of the diverse algorithms on the four datasets in
% Figure (\red{TBD}), along with the corresponding average Mean Squared Prediction Error (MSPE) values, as summarized in
Table \ref{tab:MSPE_compare}, where we emphasize \textbf{the best outcomes} using bold text and underline \underline{the second-best outcomes}.

Through a thorough analysis of the depicted visual representations and the provided tabulated data, we can deduce the following insights.

\begin{table}[h!]
    \centering
    \caption{MSPE comparison with different tensor-based algorithms on four datasets.}
    \footnotesize
    \renewcommand{\arraystretch}{1.5} % ÐÐ¼ä¾àÉèÖÃÎªÄ¬ÈÏµÄ1.5±¶
    \setlength{\tabcolsep}{2pt} % ÁÐ¿í
    \begin{tabular}{l|cccc}
    \hline & $\begin{array}{c}\text {SYN} \\
    \left(\times 10^{-2}\right)\end{array}$ & $\begin{array}{c}\text {USHCN} \\
    \left(\times 10^{-1}\right)\end{array}$ & $\begin{array}{c}\text {NASDAQ100} \\
    \left(\times 10^{-2}\right)\end{array}$ & $\begin{array}{c}\text {CCDS} \\
    \left(\times 10^{-1}\right)\end{array}$ \\
    \hline
    \hline ARIMA & $\underline{0.670}$ & 19.64 & $\underline{5.28}$ & 5.834 \\
    \hline DeepAR & 1.762 & 2.834 & 6.649 &  1.142\\
    \hline MOAR & 6.713 & 1.939 & 16.50 & 28.17 \\
    \hline MCAR & 6.152 & 1.834 & 8.12 & $\mathbf{0.562}$ \\
    \hline BHT-ARIMA & 4.161 & $\underline{1.765}$ & 8.46 & 1.567 \\
    \hline \hline LOTAP & $\mathbf{0.461}$ & $\mathbf{1.388}$ & $\mathbf{5.23}$ & $\underline{0.569}$ \\
    \hline
    \end{tabular}
    \label{tab:MSPE_compare}
\end{table}

\textbf{For synthetic dataset:} In the context of the SYN dataset, our LOTAP algorithm attains the highest level of predictive performance, with the traditional ARIMA algorithm closely following. Notably, the three temporal forecasting algorithms rooted in Tucker decomposition, namely MCAR, MOAR, and BHT-ARIMA, all exhibit subpar performance and lag considerably behind the traditional ARIMA algorithm, which refrains from employing tensor decomposition methodologies.

This discrepancy can be attributed to the fact that the data within the dataset exclusively adheres to a low-rank nature when considered under the framework of the t-SVD algorithm, whereas it does not manifest a low-rank structure within the context of Tucker decomposition. Consequently, employing a low-rank approximation through the Tucker decomposition leads to a substantial loss of pivotal information embedded within the data itself. This accounts for the diminished predictive performance compared to the ARIMA algorithm.

These findings further corroborate our earlier analysis, which underscored that time series forecasting algorithms based on t-SVD exhibit a broader spectrum of applicability compared to algorithms aligned with Tucker decomposition principles. Conversely, our proposed LOTAP algorithm leverages t-SVD to extract crucial data information and seamlessly integrates it with the AR model to identify latent temporal correlation patterns within the data. This integration contributes to its superior predictive performance surpassing that of the ARIMA algorithm.

\textbf{For real-world datasets:} Across the scope of these three real-world datasets, our LOTAP algorithm excels in predictive performance within the USHCN and NASDAQ100 datasets, while securing the second best performance in the CCDS dataset.

Significantly, a number of algorithms leveraging tensor decomposition techniques markedly outperform the non-tensor-based ARIMA and DeepAR algorithms in both the USHCN and CCDS datasets. This outcome underscores the importance of considering multicollinearity in the analysis of higher-order time series data, further affirming the value of incorporating tensor-based approaches.

Additionally, it's worth noting that the MOAR algorithm consistently falls short in predictive performance when juxtaposed with the MCAR algorithm. This divergence can be attributed to the greater regularity in terms within the MCAR algorithm, which consequently yields more stable results. This finding underscores the importance of introducing additional regularization terms to our LOTAP model.

Notably, in the NASDAQ100 dataset, the ARIMA algorithm once again secures the second best predictive performance, even surpassing algorithms like MCAR that leverage Tucker decomposition. This peculiarity is likely due to the fact that the time series data in the NASDAQ100 dataset does not fully conform to low-rank characteristics within the Tucker decomposition framework. This supposition gains further support from a low-rankness test conducted on the NASDAQ100 dataset via Tucker decomposition, revealing that around $30\%$ of the time points do not exhibit low-rank properties within this paradigm. This discovery is significant, as it underscores the existence of real-world data scenarios where established time series forecasting models rooted in Tucker decomposition struggle to provide effective solutions. In such cases, the LOTAP algorithm emerges as the solitary contender capable of accurately extracting key data information.

\begin{figure}[ht]
    \centering
    \includegraphics[width=0.75\linewidth]{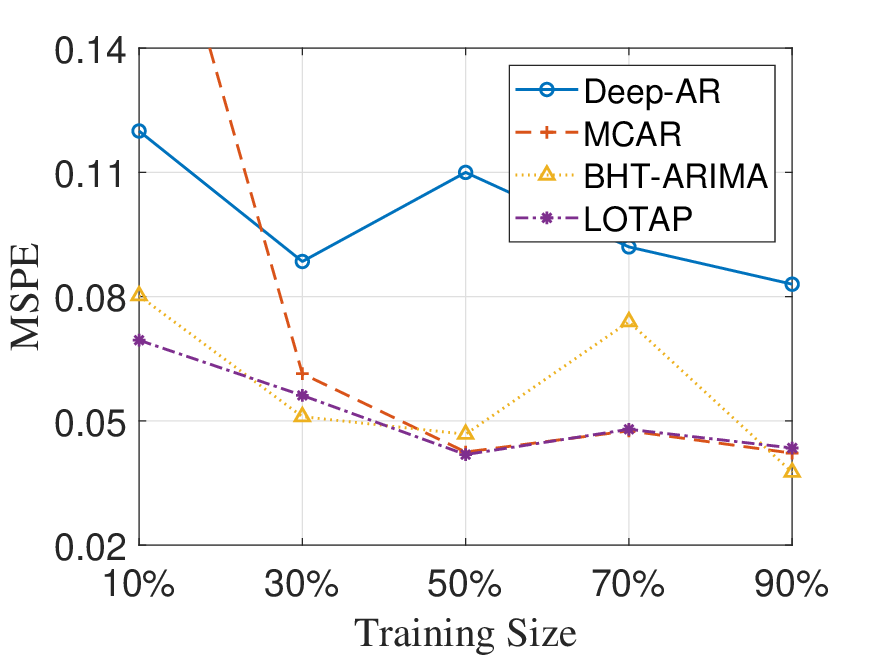}
    \caption{Comparison of algorithms accuracy under different percentage size training sets in CCDS dataset. In the case of using $10\%$ training size, the MSPE of MCAR is too large ($0.213$) so we truncate its curve there.}
    \label{fig:training size}
\end{figure}

\paragraph{4) Evaluation of training size}
To provide further validation regarding the impact of training set size on various algorithms, we conducted a series of experiments on the CCDS dataset, systematically varying the training set size from $10\%$ to $90\%$ of the tensor slices. For each dataset configuration, the remaining $10\%$ of the time slices were reserved as the test set, and the experimental outcomes are depicted in Figure \ref{fig:training size}. 
{This design not only ensures comparability across different training lengths but also avoids the unfair bias that would arise if much longer forecasting horizons were imposed.}

Importantly, due to the substantial MSPE values exhibited by the MOAR and ARIMA algorithms (both exceeding 0.5), the corresponding curves for these algorithms have been excluded from the graph. Figure 5 highlights that even with a mere $10\%$ of the training data (equivalent to 14 time points), the LOTAP algorithm demonstrates the capability to forecast the test set data with minimal error, a feat unattainable by the MCAR algorithm. Furthermore, it's worth noting that while the BHT-ARIMA algorithm marginally outperforms other algorithms in forecasting accuracy when utilizing $30\%$ and $90\%$ of the training data, its performance under different percentages of the training set reveals its inherent instability with respect to varying training set sizes.

\begin{table}[h!]
    \centering
    \caption{Time cost (in milliseconds) comparison with different tensor-based algorithms on four datasets.}
    \footnotesize
    \renewcommand{\arraystretch}{2} % ÐÐ¼ä¾àÉèÖÃÎªÄ¬ÈÏµÄ1.5±¶
    \begin{tabular}{l|cccc}
    \hline Time(ms) & \text {SYN}& \text {USHCN} & \text {NASDAQ100} & \text {CCDS} \\
    \hline
    \hline MOAR & 1587.3 & 17.4 & 297.8 & 147.6 \\
    \hline MCAR & 1324.2 & 16.2 & 301.9 & 161.7 \\
    \hline BHT-ARIMA & 24575.1 & 137.4 & 668.6 & 249.3 \\
    \hline
    \hline LOTAP & $\mathbf{264.6}$ & $\mathbf{7.9}$ & $\mathbf{58.2}$ & $\mathbf{36.2}$ \\
    \hline
    \end{tabular}
    \label{tab:Timecost_compare}
\end{table}

\paragraph{5) Time cost comparison}
We now undertake a comparative analysis of the average time dimension for forecasting in relation to three Tucker decomposition-based algorithms.
Note that data-driven methods such as DeepAR require a lot of time to train the data, so for fairness reasons we only consider these four model-driven algorithms.
The Table \ref{tab:Timecost_compare} provides a comprehensive overview of the average time expenditure per iteration during the forecasting process. Notably, our proposed LOTAP algorithm showcases minimal time overhead across all four datasets. In contrast, the BHT-ARIMA algorithm remarkably outperforms several alternative algorithms in terms of time overhead, primarily due to its multidimensional tensor strategy, which enhances tensors from three-dimensional to four-dimensional. These experimental findings align with the complexity analysis outcomes presented in Table \ref{table: complexity comparison}, thereby substantiating the assertion that the time series forecasting algorithm predicated on t-SVD yields significantly reduced time overhead compared to its Tucker decomposition-based counterpart.

\section{Conclusion}
In this paper, we have established the innovative LOTAP algorithm tailored for enhancing the realm of higher-order time series forecasting. Leveraging the truncated t-SVD technique, LOTAP extracts the intrinsic time-correlation model from data and combines it with an AR time series model for forecasting. Comparative analysis demonstrates LOTAP's superiority in applicability and forecasting speed over Tucker-based methods. Numerical experiments validate its effectiveness on diverse datasets, marking LOTAP as a promising approach for third-order time series forecasting.

\section*{Credit authorship contribution statement}
{\bf Haoning Wang}: Writing--review \& editing, Writing--original draft, Visualization, Validation, Software, Methodology, Data curation, Conceptualization.  {\bf Liping Zhang}: Writing--review \& editing, Supervision, Methodology, Funding acquisition, Formal Analysis, Conceptualization.

\section*{Declaration of competing interest}
The authors declare that they have no known competing financial
interests or personal relationships that could have appeared to
influence the work reported in this paper.

\section*{Data availability}
All codes can be found in https://github.com/whn18/LOTAP.

\section*{Acknowledgements}
This work was supported by the National Nature Science
Foundation of China (Grant No. 12571323, 12171271). %The authors thank the editor and referees for comments and suggestions that were helpful for improving this manuscript.

\section*{References}
\bibliographystyle{elsarticle-num}
\bibliography{ref}

\end{document}